\documentclass[]{interact}

\usepackage{braket,amsfonts}
\usepackage{subdepth}

\usepackage{array}

\usepackage[caption=false]{subfig}

\usepackage{algorithmic}
\usepackage{algorithm}

\usepackage{graphicx,epstopdf}

\usepackage{url}
\usepackage{amsopn}

\usepackage{graphics}
\usepackage{amssymb}
\usepackage{amsmath}
\usepackage{xcolor}
\usepackage{graphicx}
\usepackage{todonotes}
\usepackage{mathrsfs}

\usepackage{hyperref}
\hypersetup{
  breaklinks,
  colorlinks,
  citecolor={green!50!black},
  urlcolor={red!50!black},
  linkcolor={green!50!black}
}
\usepackage[nameinlink]{cleveref}
\crefname{subsection}{section}{sections}
\Crefname{algorithm}{Algorithm}{Algorithms}

\Crefname{ALC@unique}{Line}{Lines}

\usepackage[numbers,sort&compress]{natbib}
\makeatletter
\renewcommand\NAT@bibsetnum[1]{\settowidth\labelwidth{\@biblabel{#1}}%
   \setlength{\leftmargin}{\bibindent}\addtolength{\leftmargin}{\dimexpr\labelwidth+\labelsep\relax}%
   \setlength{\itemindent}{-\bibindent}%
   \setlength{\listparindent}{\itemindent}
\setlength{\itemsep}{\bibsep}\setlength{\parsep}{\z@}%
   \ifNAT@openbib
     \addtolength{\leftmargin}{\bibindent}%
     \setlength{\itemindent}{-\bibindent}%
     \setlength{\listparindent}{\itemindent}%
     \setlength{\parsep}{0pt}%
   \fi
}
\makeatother

\usepackage{etoolbox}
\makeatletter
\patchcmd{\NAT@test}{\else \NAT@nm}{\else \NAT@hyper@{\NAT@nm}}{}{}
\makeatother

\newcommand{\Kc}{{\cal K}}

\newcommand{\norm}[1]{\|#1\|}

\newtheorem{theorem}{Theorem}[section]
\crefname{theorem}{Theorem}{Theorems}
\newtheorem{cor}{Corollary}[section]
\crefname{cor}{Corollary}{Corollaries}
\newtheorem{lem}{Lemma}[section]
\crefname{lem}{Lemma}{Lemmas}
\newtheorem{rem}{Remark}[section]
\crefname{rem}{Remark}{Remarks}
\newtheorem{ex}{Example}
\crefname{ex}{Example}{Examples}

\newcommand{\pmat}[1]{\begin{pmatrix}#1\end{pmatrix}} 
\newcommand{\bmat}[1]{\begin{bmatrix}#1\end{bmatrix}} 
\newcommand{\inv}{^{-1}}
\newcommand{\T}{^T\!}

\newcommand{\calH}{{\cal H}}
\newcommand{\calS}{{\cal S}}

\usepackage{fancyhdr}

\title{A semi-conjugate gradient method for solving
       unsymmetric \\ positive definite linear systems}

\author{Na Huang%
  \thanks{Department of Applied Mathematics, College of Science,
  China Agricultural University,
    Beijing, China.
    E-mail: hna@cau.edu.cn.
    Research partially supported by National Natural Science Foundation of China (No.\,12001531).
  }
  \and Yu-Hong Dai%
  \thanks{LSEC, Academy of Mathematics and Systems Science,
    Chinese Academy of Sciences, Beijing, China.
    E-mail: dyh@lsec.cc.ac.cn.
  }
  \and Dominique Orban%
  \thanks{GERAD and
    Department of Mathematics and Industrial Engineering,
    Polytechnique Montr\'eal, QC, Canada.
    E-mail: dominique.orban@gerad.ca.
    Research partially supported by an NSERC Discovery Grant.
  }
  \and Michael A. Saunders%
  \thanks{Systems Optimization Laboratory,
    Department of Management Science and Engineering,
    Stanford University, Stanford, CA, USA.
    E-mail: saunders@stanford.edu.
  }
}

\begin{document}

\maketitle

\begin{abstract}
The conjugate gradient (CG) method is a classic Krylov subspace method for solving symmetric positive definite linear systems.
We introduce an analogous semi-conjugate gradient (SCG) method for unsymmetric positive definite linear systems. Unlike CG, SCG requires the solution of a lower triangular linear system to produce each semi-conjugate direction. We prove that SCG is theoretically equivalent to the full orthogonalization method (FOM), which is based on the Arnoldi process and converges in a finite number of steps. Because SCG's triangular system increases in size each iteration, we study a sliding window implementation (SWI) to improve efficiency, and show that the directions produced are still locally semi-conjugate. A counter-example illustrates that SWI is different from the direct incomplete orthogonalization method (DIOM), which is FOM with a sliding window. Numerical experiments from the convection-diffusion equation and other applications show that SCG is robust and that the sliding window implementation SWI allows SCG to solve large systems efficiently.
\end{abstract}


\begin{amscode}
  15A06, 65F10, 65F25, 65F50
\end{amscode}


\section{Introduction}
We consider numerical methods for solving linear systems
\begin{equation}\label{a1}
   Ax=b,
\end{equation}
where $A\in\mathbb{R}^{n\times n}$ is unsymmetric positive definite. Such a matrix $A$ is positive definite if $x\T Ax>0$ holds for all nonzero $x\in\mathbb{R}^n$ \citep{Golub1979}. This is true if and only if its symmetric part $(A+A^T)/2$ is symmetric positive definite.

Krylov subspace methods 
seek an approximate solution $x_k$ from the affine subspace $x_0+\Kc_k(A,r_0)$, where $x_0$ is an arbitrary initial point, $r_0=b-Ax_0$ is the initial residual, and $\Kc_k(A,r_0)$ is the Krylov subspace
$$\Kc_k(A,r_0)={\rm span}\{r_0,Ar_0,A^2r_0,\dots,A^{k-1}r_0\},$$
which we denote by $\Kc_k$ when there is no ambiguity.

One of the best known Krylov subspace methods is the conjugate gradient (CG) method \citep{Hestenes1952}, which was derived in 1952 to solve sparse symmetric positive definite linear systems. When combined with a suitable preconditioning, CG has many successful applications in science and engineering.
If $A$ is unsymmetric or rectangular,
CG could be applied to the normal equations $A\T Ax=A\T b$.
It is numerically preferable to apply LSQR \citep{Paige1982} or LSMR \citep{FonS2011} to $\min \norm{Ax-b}_2^2$, but the squared condition may lead to excessive iterations on compatible systems $Ax=b$.

A variety of other methods
have been developed to deal with the square unsymmetric case, such as the generalized CG-type methods~(e.g., Orthodir, Orthomin \citep{Young1980,Vinsome1976,axelsson1980conjugate,axelsson1980generalized}), the biconjugate gradient (BiCG) algorithm and its variations~(e.g., CGS, BiCGSTAB, QMR, CSBCG \citep{Sonneveld1989,Van1992,Freund1991,bank1993analysis,Tong2000}), the Manteuffel-Chebyshev iterations \citep{Manteuffel1977,Manteuffel1978,Van1981}, and other generations based on orthogonal factorizations, Lanczos process~(e.g., USYMLQ, USYMQR, LSQR, BiLQ \citep{Saunders1988,Paige1982,Montoison2020}).

Arnoldi's method \citep{arnoldi1951principle} was introduced in 1951 to deal with square unsymmetric matrices. This 
is an algorithm for constructing an orthonormal basis of the Krylov subspace $\Kc_k$. Subsequently, based on the Arnoldi process or its variations \citep{Saad}, several Krylov subspace methods for solving unsymmetric linear systems were established, such as the generalized minimum residual method (GMRES) \citep{Saad1986}, the direct quasi-GMRES method (DQGMRES) \citep{Saad1996}, the full orthogonalization method (FOM) \citep{Saad}, and the direct incomplete orthogonalization method (DIOM) \citep{Saad}.

Unlike finding an approximation from a Krylov subspace, \citet{yuan2004} sought an approximation following some semi-conjugate directions and  presented two semi-conjugate direction (SCD) methods for general linear systems. They showed that SCD has no breakdown for real positive definite systems. Later, \citet{dy2004} further studied SCD methods and introduced a new implementation for generating the semi-conjugate directions using only the latest $m$ conjugate directions, where $m$ is a given positive integer. SCD methods also have been considered for solving nonlinear systems of equations and finding saddle points of functions, which called pseudo-orthogonal direction methods in \cite{Oleg1982a,Oleg1982b}.

Here we also focus on SCD methods,
taking the first direction to be $r_0$ as in CG.
Hence, we call it the semi-conjugate gradient (SCG) method. We show that SCG is theoretically equivalent to FOM. 
SCG needs to solve a lower triangular linear system
of increasing size at each step.
To improve efficiency, we also study the sliding window implementation of SCG (SWI) and show that SWI still belongs to the set of Krylov subspace methods and will not break down. In contrast to the sliding window implementation of FOM (i.e., DIOM), the directions produced by SWI are locally  semi-conjugate. This illustrates that SWI is different from DIOM. Several numerical experiments on linear systems from the convection-diffusion equation and other applications show that SWI often solves problems
more efficiently than SCG and DIOM.

The paper is organized as follows. In \cref{sec:scg}, we introduce our semi-conjugate gradient method and provide some important properties. In \cref{subsec:scg&fom}, we prove that SCG is theoretically equivalent to FOM. The sliding window implementation of SCG and its convergence analysis are provided in \cref{sec:swi}. A counter-example in \cref{subsec:swi&diom} illustrates that SWI and DIOM are different. Numerical experiments are reported in \cref{sec:numres}. Conclusions and future work are summarized in \cref{sec:con}.

\subsection*{Notation}
For a matrix $W\in\mathbb{R}^{n\times n}$, $\calH_W=(W+W^T)/2$ and $\calS_W=(W-W^T)/2$ denote the symmetric and skew-symmetric parts of $W$. $\lambda(W)$ and $\rho(W)$ denote an arbitrary eigenvalue and the spectral radius of $W$. $\lambda_{\min}(W)$ and $\lambda_{\max}(W)$ denote the minimum and maximum eigenvalues of a symmetric matrix $W$.
A vector $e_k$ is the \(k\)th column of an identity matrix.
The solution of $Ax=b$~\eqref{a1} is denoted by $x_\star$.
The $k$th approximation to $x_\star$ is $x_k$, and the corresponding error is $d_k = x_k - x_\star$.  The 2-norm $\norm{v}$ is used for vectors $v$.

\section{The semi-conjugate gradient method}\label{sec:scg}

In this section, we introduce the semi-conjugate gradient method SCG to solve unsymmetric positive definite linear systems~\eqref{a1}. The method is summarized in \Cref{algscg}. Without loss of generality, we assume that $x_0=0$ and then $r_0=b$.

\begin{algorithm}
\caption{SCG: The semi-conjugate gradient method}
\label{algscg}
\begin{algorithmic}[1]
\STATE{Given $x_0=0$, set $k=0$, $r_0=b$, $p_0=r_0$ and $q_0=Ap_0$.}
\WHILE{a stopping condition is not satisfied}
\STATE{Compute the step size
  $\alpha_k={p_k^Tr_k}/{p_k^Tq_k}$.
\\ Update $x_{k+1}=x_k+\alpha_k p_k$ and $r_{k+1}=r_k-\alpha_k q_k.$}

\STATE{Form $v_{k+1}=Ar_{k+1}$, $L_{k+1}=P_{k+1}^TQ_{k+1}$
        and solve $L_{k+1}\lambda_{k+1} = P_{k+1}^Tv_{k+1}$,
\\ where $P_{k+1}=\bmat{p_0 & p_1 &\dots& p_k}$, $Q_{k+1}=\bmat{q_0 & q_1 &\dots& q_k}$.}
\STATE{Update
\vspace*{-15pt}
     \begin{equation}
        p_{k+1} = r_{k+1}-P_{k+1}\lambda_{k+1},  \qquad
        q_{k+1} = v_{k+1}-Q_{k+1}\lambda_{k+1}.  \label{ac2}
     \end{equation}}
\vspace*{-18pt}
\STATE{Increment $k$ by $1$.}
\ENDWHILE
\end{algorithmic}
\end{algorithm}

As shown in \cref{l1} below, $L_{k+1}$ is lower triangular. The system $L_{k+1}\lambda_{k+1} = P_{k+1}^Tv_{k+1}$ has the form
$$
  \pmat{p_{0}^Tq_0
     \\ p_{1}^Tq_0 &p_{1}^Tq_1
     \\ p_{2}^Tq_0 &p_{2}^Tq_1 &p_{2}^Tq_2
     \\ \vdots & \vdots & \vdots & \ddots
     \\ p_k^Tq_0 &p_k^Tq_1 &p_k^Tq_2 &\ \dots & p_k^Tq_k
  }
  \pmat{\lambda_{k+1}^{(1)}
     \\ \lambda_{k+1}^{(2)}
     \\ \lambda_{k+1}^{(3)}
     \\ \vdots
     \\ \lambda_{k+1}^{(k+1)}
  }
  =
  \pmat{p_{0}^{T}v_{k+1}
     \\ p_{1}^{T}v_{k+1}
     \\ p_{2}^{T}v_{k+1}
     \\ \vdots
     \\ p_k^{T}v_{k+1}
  },
$$
where $\lambda_{k+1}^{(i)}$ denotes the $i$th component of $\lambda_{k+1}$. Hence
\begin{align*}
  \lambda_{k+1}^{(1)} & = \frac{p_{0}^{T}v_{k+1}}{p_{0}^Tq_0}, \\
  \lambda_{k+1}^{(i)} & = \frac{p_{i-1}^T(v_{k+1}-\lambda_{k+1}^{(1)}q_0-\lambda_{k+1}^{(2)}q_1-\dots-       \lambda_{k+1}^{(i-1)}q_{i-2})}      {p_{i-1}^Tq_{i-1}}
  \quad (i > 1).
\end{align*}
\Cref{algpq} combines this with~\eqref{ac2} to compute directions $p_{k+1}$ and $q_{k+1}$.




\begin{algorithm}
\caption{Computation of \(p_{k+1}\) and \(q_{k+1}\)}
\label{algpq}
\begin{algorithmic}
\STATE{Assume that $p_0,\dots,p_k$, $q_0,\dots,q_k$ and $r_{0},\dots,r_{k+1}$ have been computed.}
\STATE{Set $p \leftarrow r_{k+1}$ and $q \leftarrow v_{k+1}$.}
 \FOR{$i=1, 2, \dots, k+1$}
\STATE{Compute $\lambda_{k+1}^{(i)} =  p_{i-1}^Tq/p_{i-1}^Tq_{i-1}$.}
\STATE{ Set $p\leftarrow  p - \lambda_{k+1}^{(i)}p_{i-1}$
         and $q \leftarrow  q - \lambda_{k+1}^{(i)}q_{i-1}$.}
\ENDFOR
\STATE{Set $p_{k+1} = p$ and $q_{k+1} = q$.}
\end{algorithmic}
\end{algorithm}

As stated, \Cref{algscg} and \Cref{algpq} together form a special case of \citep[Algorithm~2.3]{dy2004} in which \(p_0 = r_0\).
In other words, SCG is a special case of SCD \citep[Algorithm~4.1]{yuan2004}, though the latter does not explicitly use \(q_k\).
We now derive some further
important properties of SCG and show that it generates the same iterates as FOM. Thus, SCG theoretically follows the iterations of CG if $A$ is SPD.

\subsection{Convergence analysis of SCG}

We provide properties of SCG and prove that it converges in a finite number of steps.

\begin{lem}\label{lpq}
The sequence $\{p_k,\,q_k\}$ produced by SCG satisfies $q_k=Ap_k$.
\end{lem}

\begin{proof}
For $k=0$, we have $q_0=Ap_0$ and $Q_1=AP_1$.
For $k\ge1$, by recursion,~\eqref{ac2}, we have
\begin{align*}
q_k=v_k-Q_k\lambda_k=Ar_k-AP_k\lambda_k=Ap_k.
\end{align*}
\end{proof}

\begin{lem}\label{l1}
The matrices $L_k~(k\ge1)$ are nonsingular and lower triangular with positive diagonal elements.
\end{lem}

\begin{proof}
The proof is by induction on $k$. 
From \cref{lpq} and the fact that $A$ is positive definite, $L_1=P_1^T Q_1=p_0^TAp_0>0$, so the result holds for $k=1$. Assume $L_k$ possesses the desired property.
By definition of \(L_k\),
\begin{align*}
   P_k\T q_k =P_k\T (v_k-Q_kL_k\inv P_k\T v_k)
              =P_k\T v_k-P_k\T Q_kL_k\inv P_k\T v_k =0,
\end{align*}
and from \cref{lpq} we see that
\begin{align*}
   L_{k+1} &=P_{k+1}^TQ_{k+1} = \bmat{P_k&p_k}\T \bmat{Q_k&q_k}
            =\pmat{P_k\T Q_k & P_k\T q_k
                \\ p_k\T Q_k & p_k\T q_k}
\\         &=\pmat{L_k & 0
                \\ p_k\T Q_k & p_k\T Ap_k}
\end{align*}
is also nonsingular and lower triangular with positive diagonal elements.
\end{proof}

\begin{rem}\label{rem1}
As $A$ is positive definite, for $p_k\neq 0$ we get $p_k^Tq_k = p_k\T Ap_k > 0$. The step size $\alpha_k$ in SCG is therefore well defined. In addition, from \cref{l1}, we see that $L_{k+1}$ is nonsingular if $p_j \ne 0$, $j = 0, \dots, k$. Therefore, SCG will not break down as long as $p_k\neq 0$.
\end{rem}

From \Cref{lpq,l1} we immediately obtain the following result.

\begin{cor}\label{cor1}
   For all $j>i\ge0$, it holds that $ p_i^Tq_j = p_i^TAp_j =0$.
\end{cor}




\begin{lem}\label{l2}
After $k$ iterations of SCG we have
\begin{align}
   r_k^Tq_k &= p_k^Tq_k, \label{eq:rq=pq}
\\   P_k^Tr_k &= 0,\label{ac7}
\\ r_i^Tr_k &= 0,\quad i=0,1, \dots,k-1. \label{ac8}
\end{align}
\end{lem}

\begin{proof}
It follows from \eqref{ac2} that
\begin{equation}\label{ac6}
  r_k=P_kL_k\inv P_k\T Ar_k+p_k.
\end{equation}
This along with \Cref{cor1} leads to
\begin{eqnarray*}
   r_k^Tq_k=(P_kL_k\inv P_k\T Ar_k+p_k)^Tq_k
 = r_k\T A^TP_kL_k^{-T}P_k^Tq_k+p_k^Tq_k
 = p_k^Tq_k.\label{ac3}
\end{eqnarray*}
Therefore, \eqref{eq:rq=pq} holds.

The proof of \eqref{ac7} and \eqref{ac8} is by induction on $k$. For $k=1$ with $p_0=r_0$, 
\begin{align*}
   r_0^Tr_1=r_0^T(r_{0}-\alpha_0 q_0)
  &=r_0^Tr_{0}-\frac{ p_0^Tr_0 }{ p_0^Tq_0  } r_0^Tq_0
  =r_0^Tr_{0}-\frac{ r_0^Tr_0 }{ r_0^Tq_0  } r_0^Tq_0 = 0,
\end{align*}
and \(P_1^Tr_1 =p_0^Tr_1=r_0^Tr_1=0\).
Hence,~\eqref{ac7} and~\eqref{ac8} hold for $k=1$.

Suppose~\eqref{ac7} and~\eqref{ac8} hold for some $k\ge1$. Then if $i<k$,
\begin{equation}\label{ac5}
  r_i^Tr_{k+1}=r_i^T(r_k-\alpha_k q_k)=r_i^Tr_k-\alpha_k r_i^Tq_k=-\alpha_k r_i^Tq_k.
\end{equation}
With \eqref{ac6} and \Cref{cor1} this yields
\begin{eqnarray*}
  r_i^Tq_k=(P_{i}L_{i}\inv P_{i}^T Ar_{i}+p_{i})^Tq_k
  =r_{i}^TA^TP_{i}L_{i}^{-T}P_{i}^Tq_k+p_{i}^Tq_k = 0.
\end{eqnarray*}
Substituting into~\eqref{ac5} gives $r_i^Tr_{k+1}=0$.

If $i=k$, we check directly that
\begin{equation}\label{ac1}
  r_k^Tr_{k+1}=r_k^T(r_k-\alpha_k q_k)=r_k^Tr_k-\frac{ p_k^Tr_k }{ p_k^Tq_k  } r_k^Tq_k.
\end{equation}
Using~\eqref{ac2} and the inductive assumption, we have
\begin{eqnarray}
   p_k^Tr_k=(r_k-P_kL_k\inv P_k\T Ar_k)^Tr_k
 = r_k^Tr_k-r_k\T A^TP_kL_k^{-T} P_k^Tr_k
 = r_k^Tr_k.\label{ac4}
\end{eqnarray}
Substituting~\eqref{eq:rq=pq} and~\eqref{ac4} into~\eqref{ac1} gives
$$r_k^Tr_{k+1}=r_k^Tr_k-\frac{ r_k^Tr_k }{ p_k^Tq_k  } p_k^Tq_k=0,$$
so that~\eqref{ac8} also holds for $k+1$.

By \Cref{cor1} and the inductive assumption, we know that
\begin{align*}
   P_k^Tr_{k+1} &= P_k^T(r_k-\alpha_k q_k)=P_k^T r_k-\alpha_k P_k^Tq_k =0,
\\ p_k^Tr_{k+1} &= p_k^T(r_k-\alpha_k q_k)=p_k^Tr_k
                              -\frac{ p_k^Tr_k }{p_k^Tq_k} p_k^Tq_k=0.
\end{align*}
This implies that
$P_{k+1}^Tr_{k+1} = 
  \pmat{P_k^Tr_{k+1} \\ p_k^Tr_{k+1}} = 0$.
\end{proof}

\smallskip

It follows from~\eqref{ac2} and \cref{l2} that
\begin{equation*}
   p_k^T r_k = (r_k-P_k\lambda_k)^T r_k
             = r_k^Tr_k - \lambda_k^T P_k^T r_k = r_k^T r_k.
\end{equation*}
This implies that $\alpha_k$ in SCG can also be defined by
$ \alpha_k = {r_k^Tr_k}/{p_k^Tq_k}$, which saves the computation of \(p_k^T r_k\) and allows us to reuse the quantity \(r_k^T r_k\) that typically appears in a stopping condition.

\smallskip

\begin{lem}\label{lo1}
After $k-1$ iterations of SCG, if $p_k=0$ then $r_k=0$.
\end{lem}

\begin{proof}
Assume $p_j \ne 0,~j=0,\dots, k-1$. From the definition of $L_k$ and \Cref{lpq,l1} we know that $P_k$ has full column rank. Let
\begin{equation}\label{ao1}
  P_k = U_k\pmat{\Sigma_k \\ 0} V_k^T
\end{equation}
be the singular value decomposition (SVD), where $n\times n$ $U_k$ and $k\times k$ $V_k$ are unitary matrices, and $\Sigma_k=\mathrm{diag}\{\sigma_1,\sigma_2,\dots,\sigma_k\}$ with all $\sigma_j>0$ contains the singular values of $P_k$. It follows from \cref{l2} that
$$
0=P_k^Tr_k = V_k \pmat{\Sigma_k&0} U_k^Tr_k
           = V_k \pmat{\Sigma_k&0}
                 \pmat{\tilde{r}_k^{(1)} \\ \tilde{r}_k^{(2)}}
           = V_k \Sigma_k \tilde{r}_k^{(1)},
$$
where $U_k^Tr_k=((\tilde{r}_k^{(1)})^T,(\tilde{r}_k^{(2)})^T)^T$. This implies that $\tilde{r}_k^{(1)}=0$.

If $p_k=0$, by~\eqref{ac2} and \cref{lpq} we get $(I - P_k(P_k\T AP_k)\inv P_k\T A)r_k=0$,
i.e.,
\begin{eqnarray*}
  (U_k\T AU_k - U_k\T AP_k(P_k\T AP_k)\inv P_k\T A U_k)U_k^Tr_k=0
\end{eqnarray*}
as $U_k$ is unitary and $A$ is nonsingular. Substituting~\eqref{ao1} gives
\begin{equation}\label{ao2}\small
   \left[ U_k\T AU_k - U_k\T AU_k \pmat{\Sigma_k \\ 0}
      \left( \pmat{\Sigma_k & 0}
          U_k\T AU_k
          \pmat{\Sigma_k \\ 0}
      \right)\inv   \!\!
   \pmat{\Sigma_k & 0}
   U_k\T AU_k \right] U_k^Tr_k=0.
\end{equation}
Let
$U_k\T AU_k = \hbox{\small $\pmat{\tilde{A}_k^{(1)}& \tilde{A}_k^{(2)}
              \\ \tilde{A}_k^{(3)}& \tilde{A}_k^{(4)}}$}
           \in \mathbb{R}^{n\times n}$,
where $\tilde{A}_k^{(1)}\in\mathbb{R}^{k\times k}$.
Then~\eqref{ao2} reads 
$$0= \pmat{0 & 0
        \\ 0 & \tilde{A}_k^{(4)}
             - \tilde{A}_k^{(3)} \big( \tilde{A}_k^{(1)} \big)\inv
               \tilde{A}_k^{(2)}}
     \pmat{\tilde{r}_k^{(1)} \\ \tilde{r}_k^{(2)}}
               =\pmat{0\\
  \left[\tilde{A}_k^{(4)} - \tilde{A}_k^{(3)}\big(\tilde{A}_k^{(1)}\big)\inv \tilde{A}_k^{(2)}\right]
  \tilde{r}_k^{(2)}}.$$
Since $A$ is positive definite, we have $\tilde{r}_k^{(2)}=0$. This along with $\tilde{r}_k^{(1)}=0$ leads to $U_k^Tr_k=0$, which gives the result by the full column rank of $U_k$.
\end{proof}

\begin{rem}
Combining \Cref{rem1} with \Cref{lo1}, we see that SCG will not break down unless $r_k=0$.
\end{rem}

\smallskip

We are now ready to prove that SCG converges in a finite number of steps.

\begin{theorem}
SCG converges to the unique solution of the linear system~\eqref{a1} within $n+1$ steps if
roundoff errors are ignored.
\end{theorem}

\begin{proof}
If $r_k\neq0$ for all $k=0,1,\dots,n-1$, by \cref{l2} $r_0,\,r_1,\dots,r_{n-1}$ are orthogonal, and therefore linearly independent. Then there exist $a_0,a_1,\dots,a_{n-1}$ such that $r_{n}=\sum_{i=0}^{n-1}a_ir_i$.
\cref{l2} then yields $r_{n}^Tr_{n}=\sum_{i=0}^{n-1}a_ir_{n}^Tr_i=0$.
\end{proof}


\subsection{SCG is equivalent to FOM}\label{subsec:scg&fom}

In this section, we show that SCG and FOM are theoretically equivalent.

FOM is a Krylov subspace method introduced by \citet{Saad} in which the residual associated with the $k$th solution approximation $\hat{x}_k$ satisfies the Galerkin condition
\begin{equation}%
  \label{eq:galerkin}
  \hat{r}_k=b-A\hat{x}_k \perp \Kc_k.
\end{equation}
Given the initial guess $\hat{x}_0 = 0$ and $\beta=\|\hat{r}_0\| = \|b\|$, Arnoldi's method sets $\hat{v}_1 = \hat{r}_0/\beta$ and, for $j=1,2,\ldots,k-1$, computes
\begin{equation}\label{arnoldi}
    \left\{
\begin{array}{l}
           h_{ij}  = \hat{v}_i^TA\hat{v}_j, ~ i=1,2,\ldots,j,
\\            w_j  = A\hat{v}_j-\sum_{i=1}^j h_{ij} \hat{v}_i,
\\      h_{j+1,j}  = \|w_j\|,
\\   \hat{v}_{j+1} = w_j/h_{j+1,j}.
\end{array}
\right.
\end{equation}
This process constructs $\hat{V}_k=\bmat{\hat{v}_1\! & \!\!\hat{v}_2\!\! &\dots& \!\!\hat{v}_k}$ whose columns form an orthonormal basis of $\Kc_k$ such that
\[
    A\hat{V}_k = \hat{V}_k H_k + h_{k+1,k} \hat{v}_{k+1} e_k^T,
\]
where \(H_k\) is \(k\)\(\times\)\(k\) upper Hessenberg and \(h_{k+1,k}\) will appear in \(H_{k+1}\).
FOM seeks a solution of the form \(\hat{x}_k = \hat{V}_k y_k\).
Thus,
\[
    \hat{r}_k =
    \beta \hat{v}_1 - A \hat{V}_k y_k =
    \hat{V}_{k+1}\pmat{\beta e_1 - H_k y_k
                  \\[3pt]  - h_{k+1,k} e_k^T y_k}.
\]
By~\eqref{eq:galerkin}, the $k$th approximate solution $\hat{x}_k$ in FOM is given by
\begin{align*}
   y_k = H_k^{-1}(\beta e_1).
\end{align*}
If \(A\) is positive definite, so is each \(H_k\) in exact arithmetic because \(H_k = \hat{V}_k^T A \hat{V}_k\).
Therefore, it possesses an LU factorization without pivoting\footnote{In practice, pivoting remains advisable in general for stability.}, say $H_k=\hat{L}_k\hat{U}_k$ \citep{Golub1979} with \(\hat{L}_k\) unit lower triangular and \(\hat{U}_k\) upper triangular with positive diagonal elements, which allows us to state FOM in the form of \Cref{algfom}.


\begin{algorithm}   
\caption{FOM {\rm\citep[Algorithms \(6.4\) and \(6.8\)]{Saad}} }
\label{algfom}
\begin{algorithmic}[1]
\STATE{Given $ \hat{x}_0=0$, set $\hat{r}_0=b$, $\beta = \|\hat{r}_0\|$, and $\hat{v}_1 = \hat{r}_0/\beta$.}

\FOR{$k=1,2,\dots$}
\STATE{Compute $h_{ik}$, $i=1,2,\dots,k$ and $\hat{v}_{k+1}$ by the Arnoldi process.}

\STATE{Update the LU factorization of $H_k$, i.e., obtain the last column $u_k$ of $\hat{U}_k$.}

\STATE{Compute $\zeta_k=\{\mbox{if}~k=1~\mbox{then}~\beta, \mbox{else}~-l_{k,k-1}\zeta_{k-1}\}.$}

\STATE{Compute $\hat{p}_k = (\hat{v}_k-\sum_{i=1}^{k-1}u_{ik}\hat{p}_i) / u_{kk}.$}

\STATE{Compute $\hat{x}_k = \hat{x}_{k-1}+\zeta_k\hat{p}_k$.}

\STATE{Compute $\hat{r}_k = \hat{r}_{k-1} - \zeta_k A \hat{p}_k$.}

\ENDFOR
\end{algorithmic}
\end{algorithm}

We now state properties of FOM used to analyze its connection with SCG.

\begin{lem}\label{l3}
{\rm\citep[Proposition 6.7]{Saad}}
In FOM, 
$$
\hat{r}_k = -h_{k+1,k}e_k^Ty_k \hat{v}_{k+1}=\hat{v}_{k+1}/t_{k+1},
$$
where $y_k=H_k\inv (\beta e_1)$ and $t_{k+1}=1/(-h_{k+1,k}e_k^Ty_k)$.
\end{lem}

\begin{lem}\label{l4}
{\rm\citep[Properties on page 157]{Saad}}
In FOM,
\begin{itemize}
  \item the directions $\hat{p}_k$ are semi-conjugate, i.e.,
$\hat{p}_i^TA\hat{p}_j=0$ for $i<j$;
  \item the residual vectors $\hat{r}_k$ are orthogonal, i.e.,
  $\hat{r}_i^T\hat{r}_j = 0$ for $i\neq j$.
\end{itemize}

\end{lem}

The connection between SCG and FOM can now be summarized as follows.

\begin{theorem}\label{conne}
Assume that $\hat{r}_k$ and $\hat{p}_k$ are produced by FOM, and $r_k$ and $p_k$ are produced by SCG. Then for all $k\ge 1$, if $\hat{r}_{k-1}\neq0$ and $r_{k-1}\neq0$, there exists $a_k\neq 0$ such that
\begin{align}
   \hat{p}_k &= a_{k-1} p_{k-1},  \label{peql}
\\ \hat{r}_k &= r_k.              \label{reql}
\end{align}
\end{theorem}

\begin{proof}
We use induction on $k$. For $k=1$, as $\hat{r}_0 = r_0 = p_0$, it is easy to see that
$$
\hat{p}_{1} = u_{11}\inv \hat{v}_1 =  \frac{1}{\beta u_{11}}r_0  =  \frac{1}{\beta h_{11}}p_0.
$$
Note that $h_{11} = \hat{v}_1^TA\hat{v}_1 = (r_0^TAr_0)/\beta^2$ leads to
$$
\hat{p}_{1}=\frac{\beta}{r_0^TAr_0 }p_0 = \frac{\|r_0\|}{r_0^TAr_0 }p_0.
$$
In addition, by $\hat{r}_0 = r_0=p_0$ and $\zeta_1 = \beta$, we have
$$
\hat{r}_1 = \hat{r}_0 - \zeta_1\hat{p}_{1}
=r_0 -  \frac{\beta\|r_0\|}{r_0^TAr_0 }p_0
=r_0 -  \frac{r_0^Tr_0}{r_0^TAr_0 }p_0
=r_0 -  \frac{p_0^Tr_0}{p_0^TAp_0 }p_0
=r_0 -  \alpha_0p_0
=r_1.
$$
We then have~\eqref{peql}--\eqref{reql} satisfied for $k=1$ with $a_0 = \|r_0\|/(r_0^TAr_0)\ne 0$.

Suppose there exist constants $a_i\neq0~(0\le i\le K-2)$ such that~\eqref{peql}--\eqref{reql} are satisfied for all $1\le k < K$.  We show that~\eqref{peql}--\eqref{reql} also hold for $k = K$.
Let
$$
      \hat{P}_{K-1} = \bmat{\hat{p}_{1} & \!\!\!\dots\!\!\! & \hat{p}_{K-1}},
~\, \tilde{u}_K = (u_{1,K},\dots,u_{K-1,K})^T,
~\, D_{K-1} = {\rm diag} \{a_0,\dots,a_{K-2}\}.
$$
From the induction hypothesis, $D_{K-1}$ is nonsingular and $\hat{P}_{K-1} = P_{K-1}D_{K-1}$. With \cref{l3}, this leads to
\begin{align}
   \hat{p}_K &= u_{K,K}\inv
      \big( \hat{v}_K -
            \smash[t]{\sum\limits_{i=1}^{K-1}} u_{i,K}\hat{p}_i
      \big) \nonumber
 \\ &= u_{K,K}\inv  \big(\hat{v}_K-\hat{P}_{K-1}\tilde{u}_K\big) \nonumber
 \\ &= u_{K,K}\inv  \big(t_K\hat{r}_{K-1}-P_{K-1}D_{K-1}\tilde{u}_K\big)\nonumber
 \\ &= u_{K,K}\inv \big(t_Kr_{K-1}-P_{K-1}D_{K-1}\tilde{u}_K\big).\label{hatp}
\end{align}
From \cref{l4}, $\hat{p}_i^TA\hat{p}_K=0$ holds for all $i<K$, including $\hat{P}_{K-1}^TA\hat{p}_K=0$. Since $\hat{P}_{K-1} = P_{K-1}D_{K-1}$ and the matrix $D_{K-1}$ is nonsingular, we get
$$
P_{K-1}^TA\hat{p}_K=0.
$$
Multiplying both sides of~\eqref{hatp} by $P_{K-1}^TA$ and using SCG and \cref{lpq}, we obtain
\begin{align*}
  0  = t_KP_{K-1}^TAr_{K-1}-P_{K-1}^TAP_{K-1}D_{K-1}\tilde{u}_K
   & = t_KP_{K-1}^Tv_{K-1}-P_{K-1}^TQ_{K-1}D_{K-1}\tilde{u}_K
\\ & = t_KP_{K-1}^Tv_{K-1}-L_{K-1}D_{K-1}\tilde{u}_K,
\end{align*}
where we used the identity \(A r_{K-1} = v_{K-1}\) from \Cref{algscg}.
This shows that
$$D_{K-1}\tilde{u}_K = t_KL_{K-1}\inv P_{K-1}^Tv_{K-1} = t_K\lambda_{K-1}.$$
We substitute the above into~\eqref{hatp} and use~\eqref{ac2} to obtain
$$
  \hat{p}_K = u_{K,K}\inv \big(t_Kr_{K-1}-t_KP_{K-1}\lambda_{K-1}\big)
  =u_{K,K}\inv t_Kp_{K-1}.
$$
Therefore,~\eqref{peql} holds for $k=K$ with $a_{K-1}=u_{K,K}\inv t_K\neq0$.

It follows from FOM, the induction hypothesis and \Cref{lpq} that
\begin{equation}\label{hatr}
  \hat{r}_K = \hat{r}_{K-1} - \zeta_KA\hat{p}_K = r_{K-1} - \zeta_Ka_{K-1}Ap_{K-1} = r_{K-1} - \zeta_Ka_{K-1}q_{K-1}.
\end{equation}
This along with \cref{l4}, the induction hypothesis, and \eqref{eq:rq=pq} gives
\begin{align*}
  0 = \hat{r}_{K-1}^T\hat{r}_K = r_{K-1}^T\hat{r}_K
   &= r_{K-1}^Tr_{K-1} - \zeta_Ka_{K-1}r_{K-1}^Tq_{K-1}\\
   &= r_{K-1}^Tr_{K-1} - \zeta_Ka_{K-1}p_{K-1}^Tq_{K-1}.
\end{align*}
By~\Cref{cor1} and the positive definiteness of \(A\), \(p_{K-1}^Tq_{K-1} > 0\).
Thus,
$$
  \zeta_Ka_{K-1} = \frac{r_{K-1}^Tr_{K-1}}{p_{K-1}^Tq_{K-1}}.
$$
By~\eqref{ac2} and \Cref{l2}, we obtain
\begin{align*}
  r_{K-1}^Tr_{K-1}=(p_{K-1}+P_{K-1}\lambda_{K-1})^T r_{K-1}
   &= p_{K-1}^T r_{K-1} + \lambda_{K-1}^TP_{K-1}^T r_{K-1}
\\ &= p_{K-1}^T r_{K-1}.
\end{align*}
We then have
$
  \zeta_Ka_{K-1} = \frac{p_{K-1}^T r_{K-1}}{p_{K-1}^Tq_{K-1}} = \alpha_{K-1}.
$
Combining with~\eqref{hatr} gives
$$
\hat{r}_K = r_{K-1} - \alpha_{K-1}q_{K-1} = r_K.
$$
Hence,~\eqref{reql} also holds for $k=K$.
\end{proof}

In the following, we show that neither $\{\|x_k-x_\star\|\}$ nor $\{\|x_k\|\}$ produced by FOM or SCG is monotonic. Consider
\begin{equation*}
  A=\pmat{1 & 0 & -2
       \\ 0 & 1 & 0
       \\ 2 & 0 & 2}
     \qquad\textrm{and}\qquad
  b=\pmat{1 \\ 0 \\ 0}.
\end{equation*}
Note that \(A\) is positive definite.
With $\hat{x}_0=0$, we have $\hat{r}_0=b$ and $\beta=\|\hat{r}_0\|=1$. It follows from \eqref{arnoldi} that
$$
\begin{array}{cccc}
\hat{v}_1 = \pmat{1  \\ 0 \\ 0 },
&\hat{v}_2 = \pmat{0  \\ 0 \\ 1 },
&H_1 = 1,
&H_2 = \pmat{ 1 & -2
          \\ 2 &  2}.
\end{array}
$$
Then we have
\begin{align*}
\hat{x}_1 &= \hat{V}_1y_1
= \hat{V}_1H_1^{-1}(\beta e_1)
= \pmat{ 1 & 0 & 0}\T,
\\
\hat{x}_2 &= \hat{V}_2y_2
= \hat{V}_2H_2^{-1}(\beta e_1)
= \pmat{ \frac{1}{3} &  0 & -\frac{1}{3}}\T = x_{\star}.
\end{align*}
This implies that
$\|\hat{x}_0-x_{\star}\| = \sqrt{2}/3$,  $\|\hat{x}_1-x_{\star}\| = \sqrt{5}/3$,  $\|\hat{x}_2-x_{\star}\| = 0$ and $\|\hat{x}_0\|=0$, $\|\hat{x}_1\|=1$, $\|\hat{x}_2\|=\sqrt{2}/3$. Therefore, the sequences $\{\|\hat{x}_k-x_\star\|\}$ and $\{\|\hat{x}_k\|\}$ produced by FOM (and SCG) are not monotonic.



\section{Sliding window implementation of SCG}\label{sec:swi}

In this section, we study the sliding window implementation of SCG, which is described in \Cref{algswi}.

\begin{algorithm}[ht]   
\caption{SWI: Sliding window implementation of SCG}
\label{algswi}
\begin{algorithmic}[1]
\STATE{Given $ x_0=0$ and a nonnegative integer $m$, set $r_0=b$, $p_0=r_0$ and $q_0=Ap_0$.}
\WHILE{a stopping condition is not satisfied}

\STATE{Compute the step size $\alpha_k={r_k^Tp_k}/{p_k^Tq_k}$.
   \\ Update $x_{k+1}=x_k+\alpha_k p_k$ and $r_{k+1}=r_k-\alpha_k q_k.$}

\STATE{Form $v_{k+1} = Ar_{k+1}$,
                    $L_{k+1} = P_{k+1}^TQ_{k+1}$
      and solve     $L_{k+1}\lambda_{k+1} = P_{k+1}^T v_{k+1}$,
   \\ where $P_{k+1}=\bmat{p_{k-m_k} & \!\!\dots\!\! & p_k}$,
            $Q_{k+1}=\bmat{q_{k-m_k} & \!\!\dots\!\! & q_k}$,
      $m_k=\min\{k,m\}$.}

\STATE{Update
\vspace*{-15pt}
     \begin{equation}
        p_{k+1} = r_{k+1}-P_{k+1}\lambda_{k+1},  \qquad
        q_{k+1} = v_{k+1}-Q_{k+1}\lambda_{k+1}.  \label{a6}
     \end{equation}}
\vspace*{-18pt}
\STATE{Increment $k$ by $1$.}
\ENDWHILE
\end{algorithmic}
\end{algorithm}

\citet[Algorithm 5.1]{dy2004} proposed the limited-memory left conjugate direction method, which is theoretically equivalent to SWI, but they did not provide an analysis of the method. In the following, we derive some properties of SWI and show that it is convergent under reasonable conditions.

As in the proof of \cref{lpq}, we obtain a relation between $p_k$ and $q_k$ in SWI.

\begin{lem}\label{nl3}
The sequence $\{p_k,\,q_k\}$ produced by SWI satisfies $q_k=Ap_k$.
\end{lem}

\begin{lem}\label{l5}
The SWI matrices $L_k$ $(k\ge1)$ are nonsingular and lower triangular.
\end{lem}

\begin{proof}
If $1\le k\le m+1$, it follows from \cref{l1} that $L_k$ is nonsingular and lower triangular.
Assume that $L_k$ is nonsingular and lower triangular for some \(k \geq m+1\).
Let us now show that the same is true of \(L_{k+1}\).
Let $\tilde{P}_k=\bmat{p_{k-m} & \!\!\dots\!\! & p_{k-1}}$,
    $\tilde{Q}_k=\bmat{q_{k-m} & \!\!\dots\!\! & q_{k-1}}$ and $k_m=k-m-1$. It is easy to see that
\begin{equation}\label{a5}
 \begin{array}{cc}
   P_k=\pmat{p_{k_m}&\tilde{P}_k}, \qquad& P_{k+1}=\pmat{\tilde{P}_k&p_k},
 \\[6pt]
   Q_k=\pmat{q_{k_m}&\tilde{Q}_k}, \qquad& Q_{k+1}=\pmat{\tilde{Q}_k&q_k}.
 \end{array}
\end{equation}
As $L_k$ is nonsingular and lower triangular, we have
\begin{equation}%
  \label{eq:Lk-lower-tri}
 L_k=P_k^TQ_k=\pmat{p_{k_m}^T \\[2pt] \tilde{P}_k^T}
              \pmat{q_{k_m} & \tilde{Q}_k}
 = \pmat{p_{k_m}^Tq_{k_m} & p_{k_m}^T\tilde{Q}_k
      \\[2pt] \tilde{P}_k^Tq_{k_m} & \tilde{P}_k^T\tilde{Q}_k}
 = \pmat{p_{k_m}^Tq_{k_m} & 0
      \\[2pt] \tilde{P}_k^Tq_{k_m} & \tilde{P}_k^T\tilde{Q}_k},
\end{equation}
which implies that $\tilde{P}_k^T\tilde{Q}_k$ is also nonsingular and lower triangular.
Together,~\eqref{a5} and~\eqref{eq:Lk-lower-tri} yield
\begin{equation}%
    \label{eq:Pktilde}
    \tilde{P}_k^T Q_k L_k\inv P_k^T =
    \pmat{\tilde{P}_k^T q_{k_m} & \tilde{P}_k^T \tilde{Q}_k}
    \pmat{p_{k_m}^Tq_{k_m} & 0
      \\[2pt] \tilde{P}_k^Tq_{k_m} & \tilde{P}_k^T\tilde{Q}_k}^{-1}
    \pmat{p_{k_m}^T
      \\[2pt] \tilde{P}_k^T} =
    \tilde{P}_k^T.
\end{equation}
Combining~\eqref{a6} with~\eqref{eq:Pktilde} gives
\(
    \tilde{P}_k^T q_k =
    \tilde{P}_k^T v_k - \tilde{P}_k^T Q_k L_k\inv P_k^T v_k =
    0
\).
Thus,
\begin{equation*}
  L_{k+1}=\pmat{\tilde{P}_k^T \\ p_k^T}
          \pmat{\tilde{Q}_k& q_k}
         =\pmat{\tilde{P}_k^T\tilde{Q}_k &\tilde{P}_k^T q_k
             \\ p_k^T\tilde{Q}_k&p_k^Tq_k}
         =\pmat{\tilde{P}_k^T\tilde{Q}_k&0
             \\ p_k^T\tilde{Q}_k&p_k^Tq_k}
\end{equation*}
is also lower triangular. Using \cref{nl3} and the fact that $A$ is positive definite, we have $p_k^Tq_k = p_k\T Ap_k>0$. Therefore, $L_{k+1}$ is nonsingular.
\end{proof}


With all $L_k$ nonsingular, SWI is well defined.
\cref{l5} also implies the following.

\begin{cor}\label{cor5}
For all $i \in [\max\{0,\,k-m\}, k-1]$, it holds that $p_i^Tq_k=0$.
\end{cor}

\begin{lem}\label{l6}
After $k$ iterations in SWI we have $P_k^Tr_k=0$.
\end{lem}

\begin{proof}
If $k\le m+1$, it follows from \cref{l2} that $P_k^Tr_k=0$. Now we prove that $P_k^Tr_k=0$ also holds for $k>m+1$. The proof is by induction on $k$. Assume that $P_k^Tr_k = 0$ for some $k\ge m+1$. For the case of $k+1$, it follows from~\eqref{a5} that
$$ P_k^Tr_k = \pmat{p_{k_m}^T    \\[2pt] \tilde{P}_k^T} r_k
            = \pmat{p_{k_m}^Tr_k \\[2pt] \tilde{P}_k^Tr_k} = 0.$$
This together with \cref{cor5} yields
$$\tilde{P}_k^Tr_{k+1} = \tilde{P}_k^T(r_k-\alpha_k q_k) =
  \tilde{P}_k^Tr_k - \alpha_k \tilde{P}_k^Tq_k=0.$$
From the expression for $\alpha_k $, we have
$$p_k^Tr_{k+1} = p_k^T(r_k-\alpha_k q_k) = p_k^Tr_k - \alpha_k p_k^Tq_k
  = p_k^Tr_k - \frac{r_k^Tp_k}{p_k^Tq_k} p_k^Tq_k = 0.$$
This shows that
$$
  P_{k+1}^Tr_{k+1} = \pmat{\tilde{P}_k^T \\ p_k^T} r_{k+1}
                   = \pmat{\tilde{P}_k^Tr_{k+1} \\ p_k^Tr_{k+1}} = 0.
$$
The proof follows by induction.
\end{proof}


\begin{rem}\label{rem6a}
It follows from \cref{l6} that
\begin{equation}\label{ao25}
  r_k^Tp_k=r_k^T(r_k-P_k\lambda_k)=r_k^Tr_k-r_k^TP_k\lambda_k=r_k^Tr_k.
\end{equation}
Hence, the step size $\alpha_k$ in SWI can also be updated by
$\alpha_k=r_k^Tr_k / p_k^Tq_k$.
\end{rem}

\begin{rem}\label{rem6b}
In the same way as for \cref{lo1}, we can prove that if $p_k=0$, then $r_k=0$. Hence, SWI will not break down unless $r_k=0$.
\end{rem}

In the following, we show that SWI is a Krylov subspace method.
\begin{lem}\label{swikrylov}
The sequence $\{x_k, r_k,\,p_k\}$ produced by SWI satisfies $x_k\in \Kc_k$ and $r_k,p_k\in\Kc_{k+1}$.
\end{lem}

\begin{proof}
For $k=0$, the results hold naturally. Assume that the results hold for some $k\ge0$. Then for $k+1$, it follows from SWI, \cref{nl3} and the induction hypothesis that
\begin{align*}
   x_{k+1} &= x_k +\alpha_k p_k \in \Kc_{k+1},
\\       r_{k+1} &= r_k-\alpha_k q_k = r_k-\alpha_k Ap_k \in \Kc_{k+2},
\\       p_{k+1} &= r_{k+1}-P_{k+1}\lambda_{k+1} \in \Kc_{k+2}.
\end{align*}
The proof follows by induction.
\end{proof}

We are now ready to establish the convergence theorem for SWI. For any $k\ge m$, let the SVD of $n \times (m+1)$ matrix $P_k$ be
\begin{equation}\label{ao5}
  P_k = U_k \pmat{\Sigma_k \\ 0} V_k^T,
\end{equation}
where
$\Sigma_k=\mathrm{diag}\{\sigma_1,\sigma_2,\dots,\sigma_{m+1}\}$ and $\sigma_j>0$.  Also let
\begin{equation}\label{ao7}
   U_k\T AU_k = \pmat{\tilde{A}_k^{(1)} & \tilde{A}_k^{(2)}
                   \\ \tilde{A}_k^{(3)} & \tilde{A}_k^{(4)}}
   \in\mathbb{R}^{n\times n}
   \quad \mbox{and} \quad
   U_k^Tr_k = \pmat{\tilde{r}_k^{(1)}
                 \\ \tilde{r}_k^{(2)}}
   \in \mathbb{R}^n,
\end{equation}
where $\tilde{A}_k^{(1)}\in\mathbb{R}^{(m+1)\times (m+1)}$ and $\tilde{r}_k^{(1)}\in\mathbb{R}^{m+1}$. As $A$ is positive definite, so is $U_k\T AU_k$. From \citep[Theorem 3.9]{axelsson1996iterative}, it follows that $\tilde{A}_k^{(1)}$ is also positive definite. We can then define the Schur complement of 
$\tilde{A}_k^{(1)}$:
\begin{equation}\label{schur}
  S_k=\tilde{A}_k^{(4)} - \tilde{A}_k^{(3)}\big(\tilde{A}_k^{(1)}\big)\inv \tilde{A}_k^{(2)}.
\end{equation}
By \citep[Theorem 3.9]{axelsson1996iterative}, we know that $S_k$ is positive definite and
\begin{equation}\small 
  (U_k\T AU_k)\inv =
  \pmat{\big(\tilde{A}_k^{(1)}\big)\inv  +
        \big(\tilde{A}_k^{(1)}\big)\inv  \tilde{A}_k^{(2)}S_k\inv
             \tilde{A}_k^{(3)}
        \big(\tilde{A}_k^{(1)}\big)\inv  &
       -\big(\tilde{A}_k^{(1)}\big)\inv  \tilde{A}_k^{(2)}S_k\inv
\\ -S_k\inv \tilde{A}_k^{(3)} \big(\tilde{A}_k^{(1)}\big)\inv  & S_k\inv }.
   \label{inver}
\end{equation}
It follows from the Courant-Fischer min-max theorem that
\begin{align*}
   \lambda \big(\calH_{S_k\inv }\big) &\le
   \lambda_{\max} \big(\calH_{S_k\inv }\big) \le
   \lambda_{\max} \big(\calH_{(U_k\T AU_k)\inv }\big)
 = \lambda_{\max} \big(\calH_{A\inv }\big),
\\
   \lambda \big(\rm{i}\calS_{S_k\inv }\big) &\le
   \lambda_{\max}\big(\rm{i}\calS_{S_k\inv }\big) \le \lambda_{\max}\big(\rm{i}\calS_{(U_k\T AU_k)\inv }\big)
 = \lambda_{\max}\big(\rm{i}\calS_{A\inv }\big).
\end{align*}
Similarly, we have
$$
   \lambda\big(\calH_{S_k\inv }\big) \ge
   \lambda_{\min}\big(\calH_{A\inv }\big)
   \qquad \mbox{and} \qquad
   \lambda\big(\rm{i}\calS_{S_k\inv }\big) \ge \lambda_{\min}\big(\rm{i}\calS_{A\inv }\big).
$$
Summing up, we have the following results.

\begin{lem}\label{min_max}
The eigenvalues of $\calH_{S_k\inv }$ and $\calS_{S_k\inv }$ satisfy the inequalities
\begin{eqnarray*}
  \lambda_{\min}     \big(\calH_{A\inv }\big)
  \le \lambda        \big(\calH_{S_k\inv }\big)
  \le \lambda_{\max} \big(\calH_{A\inv }\big)
  &\quad \mbox{and} \quad&
  \big|\lambda\big(\calS_{S_k\inv }\big)\big|
  \le \rho\big(\calS_{A\inv }\big).
\end{eqnarray*}
\end{lem}

\begin{theorem}\label{main_convergence}
If $A\in\mathbb{R}^{n\times n}$ is positive definite and $\lambda_{\min}(\calH_{A\inv })>\rho(\calS_{A\inv })$,
the sequence $\{x_k\}$ produced by SWI converges to the unique solution $x_\star$ of $Ax=b$~\eqref{a1}.
\end{theorem}

\begin{proof}
Let $d_k = x_\star-x_k$ be the error vector. Without loss of generality, we assume that $k>m$. Then $m_k=\min\{k,m\}=m$. From \cref{l6} and~\eqref{ao5}--\eqref{ao7}, we have
$$
0 = P_k^Tr_k = V_k \pmat{\Sigma_k & 0} U_k^Tr_k
             = V_k \pmat{\Sigma_k & 0}
                   \pmat{\tilde{r}_k^{(1)} \\ \tilde{r}_k^{(2)}}
             = V_k \Sigma_k \tilde{r}_k^{(1)},
$$
which leads to $\tilde{r}_k^{(1)}=0$.

It follows from SWI, \cref{nl3},~\eqref{ao25}, and $Ad_k=r_k$ that
\begin{eqnarray}
&&d_{k+1}^TAd_{k+1}
= (x_\star - x_{k+1})^TA(x_\star - x_{k+1})
= (x_\star - x_k-\alpha_k p_k)^TA(x_\star - x_k-\alpha_k p_k)\nonumber
\\[6pt]
&&= (d_k-\alpha_k p_k)^TA(d_k-\alpha_k p_k)
  = d_k\T Ad_k-\alpha_kd_k\T Ap_k-\alpha_kp_k\T Ad_k+\alpha_k^2p_k\T Ap_k\nonumber
\\[6pt]
&&=d_k\T Ad_k - \frac{r_k^Tp_k}{p_k^Tq_k}d_k\T Ap_k
             - \frac{r_k^Tp_k}{p_k^Tq_k}p_k\T Ad_k
             + \frac{(r_k^Tp_k)^2}{(p_k^Tq_k)^2}p_k\T Ap_k\nonumber
\\[6pt]
&&=d_k\T Ad_k - \frac{r_k^Tp_k}{p_k\T Ap_k}d_k\T Ap_k
             - \frac{r_k^Tp_k}{p_k\T Ap_k}p_k^Tr_k
             + \frac{(r_k^Tp_k)^2}{(p_k\T Ap_k)^2}p_k\T Ap_k\nonumber
\\[6pt]
&&=d_k\T Ad_k - \frac{r_k^Tp_k}{p_k\T Ap_k}d_k\T Ap_k
  = \left(1-\dfrac{r_k^Tp_k}{d_k\T Ad_k}\,\dfrac{d_k\T Ap_k}{p_k\T Ap_k}\right)
    \,d_k\T Ad_k \nonumber
\\[6pt] &&
  = \left(1 - \dfrac{r_k^Tr_k}{r_k\T A\inv r_k}
            \,\dfrac{d_k\T Ap_k}{p_k\T Ap_k}\right)
            \,d_k\T Ad_k.\label{recur1}
\end{eqnarray}
From~\eqref{ao5} and~\eqref{ao7}, we have
\begin{align}
   &\hspace{-20pt} U_k\T AU_k  - U_k\T AP_k(P_k\T AP_k)\inv P_k\T AU_k\nonumber
\\ &=U_k\T AU_k - U_k\T AU_k\pmat{\Sigma_k \\ 0}
   \Big(\pmat{\Sigma_k & 0} U_k\T AU_k\pmat{\Sigma_k \\ 0}\Big)\inv
   \pmat{\Sigma_k&0} U_k\T AU_k\nonumber
\\ &= \pmat{0 & 0 \\ 0 & S_k}. \label{ao10}
\end{align}
Note that $U_k$ is unitary. By SWI, \cref{nl3},~\eqref{inver}, and~\eqref{ao10}, we obtain
\begin{align*}
U_k\T Ap_k &= U_k\T A\big( I- P_k(P_k\T AP_k)\inv P_k\T A\big) r_k
\\         &= \big( U_k\T AU_k-U_k\T AP_k(P_k\T AP_k)\inv P_k\T   AU_k\big)U_k^Tr_k
\\         &= \pmat{0 & 0 \\ 0 & S_k}
              \pmat{\tilde{r}_k^{(1)} \\ \tilde{r}_k^{(2)}}
            = \pmat{0 \\ S_k\tilde{r}_k^{(2)}}
\\ \text{and \ }
  U_k^Tp_k &= U_k\T A\inv U_kU_k\T Ap_k = (U_k\T AU_k)\inv U_k\T Ap_k
  = \pmat{-\big(\tilde{A}_k^{(1)}\big)\inv   \tilde{A}_k^{(2)}\tilde{r}_k^{(2)}
        \\ \tilde{r}_k^{(2)}}.
\end{align*}
This together with~\eqref{inver} and $\tilde{r}_k^{(1)}=0$ yields
\begin{align}
   p_k\T Ap_k &= (U_k^Tp_k)^TU_k\T Ap_k =
   \pmat{-\big(\tilde{r}_k^{(2)}\big)^T\big(\tilde{A}_k^{(2)}\big)^T
                        \big(\tilde{A}_k^{(1)}\big)^{-T}&  \big(\tilde{r}_k^{(2)}\big)^T}
   \pmat{0 \\ S_k\tilde{r}_k^{(2)}} \nonumber
\\ &= \big(\tilde{r}_k^{(2)}\big)^T S_k\tilde{r}_k^{(2)},\label{ao20}
\\ d_k\T Ap_k &= r_k\T A^{-T}Ap_k =r_k^TU_kU_k\T A^{-T}U_kU_k\T Ap_k  \nonumber
\\ &= (U_k^Tr_k)^T(U_k\T AU_k)^{-T}U_k\T Ap_k
       =\big(\tilde{r}_k^{(2)}\big)^TS_k^{-T}S_k\tilde{r}_k^{(2)}.\label{ao21}
\end{align}
Substituting~\eqref{ao20}--\eqref{ao21} into~\eqref{recur1} gives
\begin{equation}\label{recur2}
 d_{k+1}^TAd_{k+1}=
  \left(1-\dfrac{r_k^Tr_k}{r_k\T A\inv r_k}\dfrac{\big(\tilde{r}_k^{(2)}\big)^TS_k^{-T}S_k\tilde{r}_k^{(2)}}
  {\big(\tilde{r}_k^{(2)}\big)^T S_k\tilde{r}_k^{(2)}}\right)\,d_k\T Ad_k.
\end{equation}

Let $v=S_k\tilde{r}_k^{(2)}\in\mathbb{R}^{n-m-1}$. It follows from \cref{min_max} and $\lambda_{\min}(\calH_{A\inv })>\rho(\calS_{A\inv })$ that
\begin{eqnarray*}
&&\dfrac{\big(\tilde{r}_k^{(2)}\big)^TS_k^{-T}S_k\tilde{r}_k^{(2)}}{\big(\tilde{r}_k^{(2)}\big)^T S_k\tilde{r}_k^{(2)}}
=\dfrac{v^T(S_k^{-T})^2v}{v^T S_k^{-T}v}
=\dfrac{v^T(S_k\inv )^2v}{v^T S_k\inv v}
=\dfrac{v^T\big(\calH_{S_k\inv }+\calS_{S_k\inv }\big)^2v}{v^T \calH_{S_k\inv }v}\\
&&=\dfrac{v^T\calH_{S_k\inv }^2v+v^T\calS_{S_k\inv }^2v}{v^T \calH_{S_k\inv }v}
\ge\lambda_{\min}(\calH_{S_k\inv })-
\dfrac{\max\Big\{\big|\lambda_{\min}(\calS_{S_k\inv }^2)\big|,\big|\lambda_{\max}(\calS_{S_k\inv }^2)\big|\Big\}}{\lambda_{\min}(\calH_{S_k\inv })}\\
&&\ge\lambda_{\min}(\calH_{A\inv })-
\dfrac{\rho(\calS_{A\inv })^2}{\lambda_{\min}(\calH_{A\inv })}>0.
\end{eqnarray*}
This along with the fact that
\begin{eqnarray*}
\dfrac{r_k^Tr_k}{r_k\T A\inv r_k}=\dfrac{r_k^Tr_k}{r_k^T\calH_{A\inv }r_k}
\ge \frac{1}{\lambda_{\max}\big(\calH_{A\inv }\big)}>0
\end{eqnarray*}
leads to
$$
 d_{k+1}^TAd_{k+1}\le
  \left(1-\frac{\big(\lambda_{\min}(\calH_{A\inv })\big)^2-\rho(\calS_{A\inv })^2}{\lambda_{\min}(\calH_{A\inv })\lambda_{\max}(\calH_{A\inv })}\right)d_k\T Ad_k.
$$
Hence when $\lambda_{\min}(\calH_{A\inv })>\rho(\calS_{A\inv })$, $d_k\T Ad_k\rightarrow0$ as $k\rightarrow\infty$. Since $A$ is positive definite, we have $d_k\rightarrow0$.
\end{proof}

If $A$ is symmetric positive definite, $\calS_{A\inv }=0$ and $\lambda_{\min}(\calH_{A\inv })>\rho(\calS_{A\inv })$ holds naturally, and SWI converges unconditionally.

If $A$ is a normal matrix, we have $A=X^*\Lambda X$, where $X$ is a unitary matrix and $\Lambda={\rm diag}\{a_1+{\rm i}b_1,\dots, a_n+{\rm i}b_n\}$ with $0<a_1 \le a_2 \le \dots \le a_n$ is the diagonal matrix containing the eigenvalues of $A$. Then 
\begin{align*}
   \lambda_{\min}(\calH_{A\inv })
   &=\frac{1}{2}\min\limits_j \left\{\frac{1}{a_j + {\rm i}b_j}
                                  + \frac{1}{a_j - {\rm i}b_j}\right\}
   =\min\limits_j\left\{\frac{a_j}{a_j^2+b_j^2}\right\},
\\ \rho(\calS_{A\inv })
   &=\frac{1}{2}\max\limits_j \left|\frac{1}{a_j + {\rm i}b_j}
                                  - \frac{1}{a_j - {\rm i}b_j}\right|
    =\max\limits_j\left\{\frac{|b_j|}{a_j^2+b_j^2}\right\}.
\end{align*}
If $|b_j|\ll a_j$ for all $1\le j \le n$, we have
\begin{eqnarray*}
\lambda_{\min}(\calH_{A\inv })\approx \frac{1}{a_n}
&\qquad \mbox{and} \qquad&
 \rho(\calS_{A\inv })\approx \max\limits_j\left\{\frac{|b_j|}{a_j^2}\right\}
 \leq \frac{1}{a_1}\max\limits_j\left\{\frac{|b_j|}{a_j}\right\}.
\end{eqnarray*}
Then from \cref{main_convergence} we know that SWI is convergent provided
$$
\frac{a_1}{a_n} > \max\limits_j\left\{\frac{|b_j|}{a_j}\right\}.
$$

\subsection{SWI is not equivalent to DIOM}\label{subsec:swi&diom}

Although SCG and FOM are equivalent, their sliding window implementations are different. DIOM, the sliding window implementation of FOM, is stated as \Cref{algdiom}.

\begin{algorithm}   
\caption{DIOM {\rm\citep[Algorithm 6.8]{Saad}} }
\label{algdiom}
\begin{algorithmic}[1]
\STATE{Given $\hat{x}_0=0$ and a positive integer $m$, set $\hat{r}_0=b$, $\beta = \|\hat{r}_0\|$, $\hat{v}_1 = \hat{r}_0/\beta$.}

\FOR{$k=1,2,\dots$}
\STATE{Compute $h_{ik}$, $i=\max\{1,k-m+1\},2,\dots,k+1$ and $\hat{v}_{k+1}$ using the
            incomplete orthogonalization process
            {\rm\citep[Algorithm 6.6]{Saad}}. }

\STATE{Update the LU factorization of $H_k$, i.e., obtain the last column $u_k$ of $\hat{U}_k$.}

\STATE{Compute $\zeta_k=\{\mbox{if}~k=1~\mbox{then}~\beta, \mbox{else}~-l_{k,k-1}\zeta_{k-1}\}$.}

\STATE{Compute $\hat{p}_k =  (\hat{v}_k -
               \sum_{i=k-m+1}^{k-1} u_{ik}\hat{p}_i) / u_{kk}$.}

\STATE{Compute $\hat{x}_k = \hat{x}_{k-1}+\zeta_k\hat{p}_k$.}

\STATE{Compute $\hat{r}_k = \hat{r}_{k-1} - \zeta_k A \hat{p}_k$.}

\ENDFOR
\end{algorithmic}
\end{algorithm}

From \Cref{algswi,algdiom} we see that the main difference between them is that in  DIOM, \citet{Saad} applies the sliding window idea to the Arnoldi vectors, whereas in SWI, we apply it to the transformed search directions $\{p_k\}$. As the following simple example shows, the directions $\{\hat{p}_k\}$ produced by DIOM do not satisfy
\begin{equation}\label{diompro}
  \hat{p}_i^TA\hat{p}_k=0~\quad~(\max\{0,\,k-m\}\le i\le k-1),
\end{equation}
which establishes that SWI and DIOM are different.

Consider
\begin{equation}
  A=\pmat{1 & 0 & 0 & 0 &-1
       \\ 0 & 1 & 0 &-1 & 0
       \\ 0 & 0 & 1 & 0 & 0
       \\ 0 & 1 & 0 & 1 & 0
       \\ 1 & 0 & 0 & 0 &2}
     \qquad\textrm{and}\qquad
  b=\pmat{1 \\ 1 \\ 1 \\ 0 \\ 0}. \label{couex1}
\end{equation}
Note that $A$ is positive definite. With $m=2$ and $\hat{x}_0=0$, we have $\hat{r}_0=b$ and $\beta=\|\hat{r}_0\|=\sqrt{3}$. It follows from DIOM that
\begin{align*}
H_4 = \pmat{ 1 & -\sqrt{\frac{2}{3}} & 0 & 0
    \\[6pt] \sqrt{\frac{2}{3}} & \frac{3}{2}         & -\frac{1}{2\sqrt{21}}&0
    \\[6pt]  0 & \frac{\sqrt{21}}{6} & \frac{17}{14} & -\frac{9}{7\sqrt{6}}
    \\[6pt]  0 & 0 & \frac{2\sqrt{6}}{7}&\frac{29}{28}},
    \quad
\hat{V}_4 = \pmat{\frac{1}{\sqrt{3}} & 0 & -\frac{2}{\sqrt{42}}
                                   & -\frac{3}{2\sqrt{7}}
    \\[6pt] \frac{1}{\sqrt{3}} & 0 & -\frac{2}{\sqrt{42}}
                                   &  \frac{2}{\sqrt{7}}
    \\[6pt] \frac{1}{\sqrt{3}} & 0 &  \frac{4}{\sqrt{42}}
                                   & -\frac{1}{2\sqrt{7}}
    \\[6pt]        0               &  \frac{1}{\sqrt{2}}
        & -\frac{3}{\sqrt{42}}     & -\frac{1}{2\sqrt{7}}
    \\[6pt]        0               &  \frac{1}{\sqrt{2}}
        & \frac{3}{\sqrt{42}}      &  \frac{1}{2\sqrt{7}}
}.
\end{align*}
The LU factors of $H_4$ are
\begin{align*}
\hat{L}_4= \pmat{1                      &   &   &
   \\[6pt] \sqrt{\frac{2}{3}}     & 1 &   &
   \\[6pt] 0 & \frac{\sqrt{21}}{13}   & 1 &
   \\[6pt] 0 & 0 &\frac{26\sqrt{6}}{114}  & 1},
   \quad
\hat{U}_4 = \pmat{1 & -\sqrt{\frac{2}{3}} & 0 & 0
   \\[6pt]    & \frac{13}{6} & -\frac{1}{2\sqrt{21}} & 0
   \\[6pt]    &              &  \frac{114}{91} & -\frac{9}{7\sqrt{6}}
   \\[6pt]    &              &                 &\frac{101}{76}}.
\end{align*}
Note that $\hat{P}_4=[\hat{p}_1,\hat{p}_2,\hat{p}_3,\hat{p}_4]=\hat{V}_4\hat{U}_4\inv $ \citep[p155]{Saad}.
Thus,
\begin{equation*}
  \begin{array}{llll}
    \hat{p}_1=\frac{1}{\sqrt{3}}\pmat{1 \\ 1 \\ 1 \\ 0 \\ 0},
    &
  \hat{p}_2=\frac{2\sqrt{2}}{13}\pmat{1 \\ 1 \\ 1 \\ 1 \\ 1},
    &
  \hat{p}_3=\frac{\sqrt{42}}{114}\pmat{-4 \\ -4 \\ 9 \\ -6 \\ 7},
  &
  \hat{p}_4=\frac{2}{101\sqrt{7}}\pmat{-69 \\ 64 \\ 8 \\ -37 \\ 40}.
  \end{array}
\end{equation*}
Since $\hat{p}_3^TA\hat{p}_4 =\frac{\sqrt{42}}{114}\cdot\frac{2}{101\sqrt{7}}\cdot19
= \frac{\sqrt{6}}{303}\neq0,$ we know that DIOM does not possess the properties in~\eqref{diompro}.

In addition, the SWI residuals do not satisfy
$$r_{i}^Tr_k=0 ~\quad~(\max\{0,\,k-m\}\le i\le k-1),$$
a property that the DIOM residuals possess \citep[p157]{Saad}. Indeed, for the same linear system~\eqref{couex1} and the same setting $m=2$ and $x_0=0$, the SWI residuals $\{r_k\}$ are
\begin{equation*}
  \begin{array}{llll}
r_1=\pmat{0\\0\\0\\-1\\-1},
~
r_2=\frac{1}{13}\pmat{-2\\ -2\\ 4\\ -3\\ 3},
~
r_3=\frac{1}{19}\pmat{3\\ -4\\ 1\\ 1\\ -1},
~
r_4=\frac{1}{15}\pmat{-1\\ 0 \\ 1 \\ 1\\ -1},
~
r_5=\frac{1}{289}\pmat{-1\\ 2\\ -1\\ 13 \\13}.
  \end{array}
\end{equation*}
This implies that $r_3^Tr_5 = -\frac{12}{5491} \neq 0$.

\section{Numerical experiments}\label{sec:numres}

We report numerical experience with SCG (\Cref{algscg}), SWI (\Cref{algswi}), FOM~(\Cref{algfom}), and DIOM~(\Cref{algdiom}).
For completeness, we include results obtained with GMRES \citep{Saad1986}, DQGMRES \citep{Saad1996}, and BICGSTAB \citep{Van1992}. All experiments were run using MATLAB R2015b on a PC with an Intel(R) Core(TM) i7-8550U CPU @~1.8GHz and 16GB of RAM.

In our implementation, $x_0 = 0$. Each method is terminated when either the number of iterations exceeds $10^4$ or
$$
{\rm Res}:=\frac{\|r_k\|}{\|r_0\|} < 10^{-6}.
$$
We compare the performance by reporting the number of iterations, the CPU time and the relative residual, denoted by ``Iter", ``CPU"  and  ``Res", respectively. For SWI, DIOM, and DQGMRES, we tested several values of the memory $m$, and denote the corresponding algorithms SWI($m$), DIOM($m$), and DQGMRES($m$), respectively.

\begin{ex}\label{ex_condiff}
\citep[Example~\(3.1.1\)]{elman-silvester-wathen-2006}
We consider the 2D convection-diffusion equation
\begin{equation*}
 -\epsilon\nabla^2u+\vec{w}\cdot\nabla u=0\quad \text{in } (-1, \, 1) \times (-1, \, 1),
\end{equation*}
with boundary conditions
$$
u(x,-1)=x,\quad
u(x,1)=0,\quad
u(-1,y) \approx -1,\quad
u(1,y) \approx 1.
$$
If $\vec{w}=(0,1)$, an exact solution is
$$
u(x,y) = x (1 - e^{\frac{y-1}{\epsilon}}) / (1 - e^{-\frac{2}{\epsilon}}),
$$
which satisfies the boundary conditions, save for the last two near $y = 1$.
\end{ex}

\begin{figure}[t]   
	\centering
	\subfloat{
		\includegraphics[width=.48\linewidth]{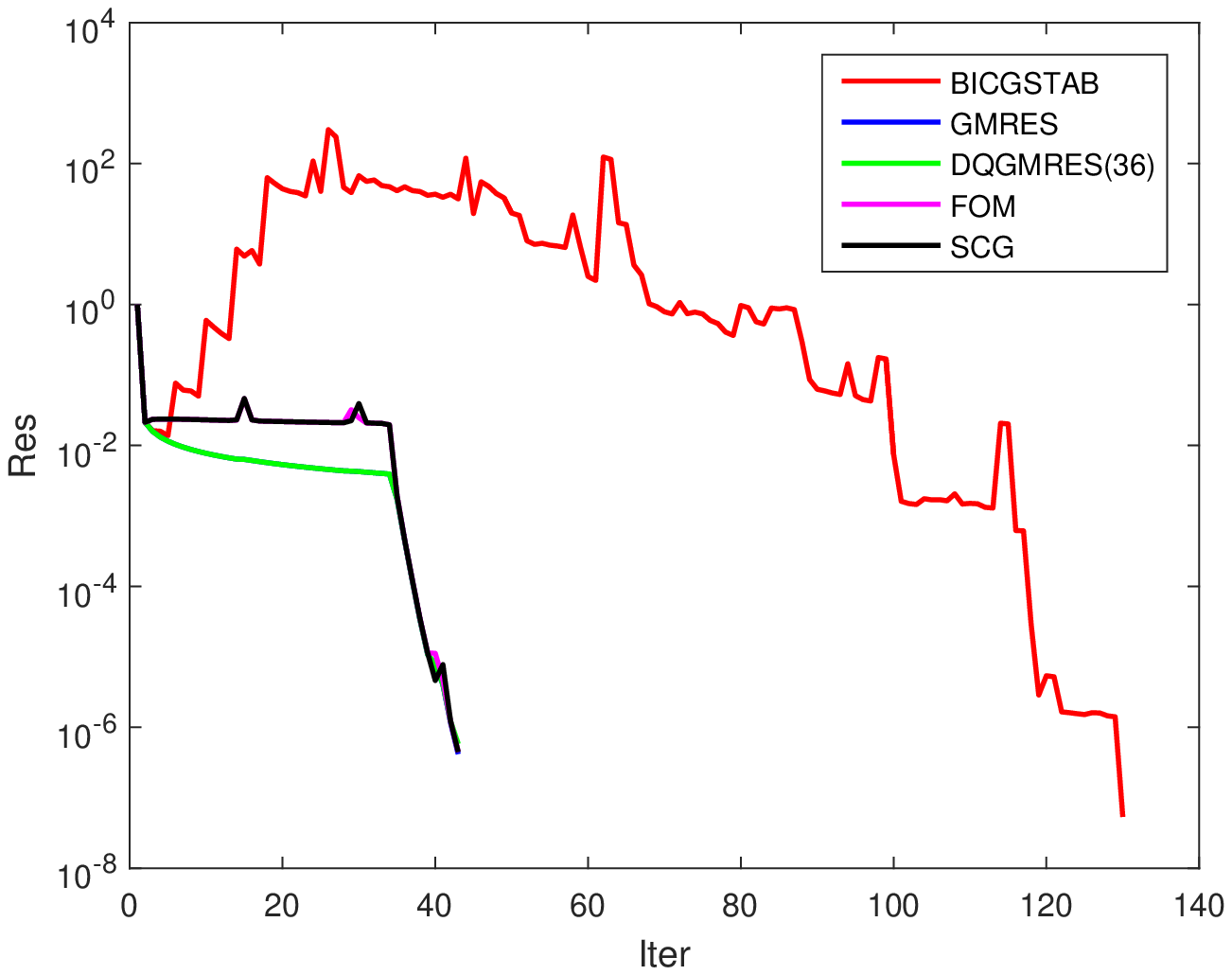}
	}
	\hfill
	\subfloat{
		\includegraphics[width=.48\linewidth]{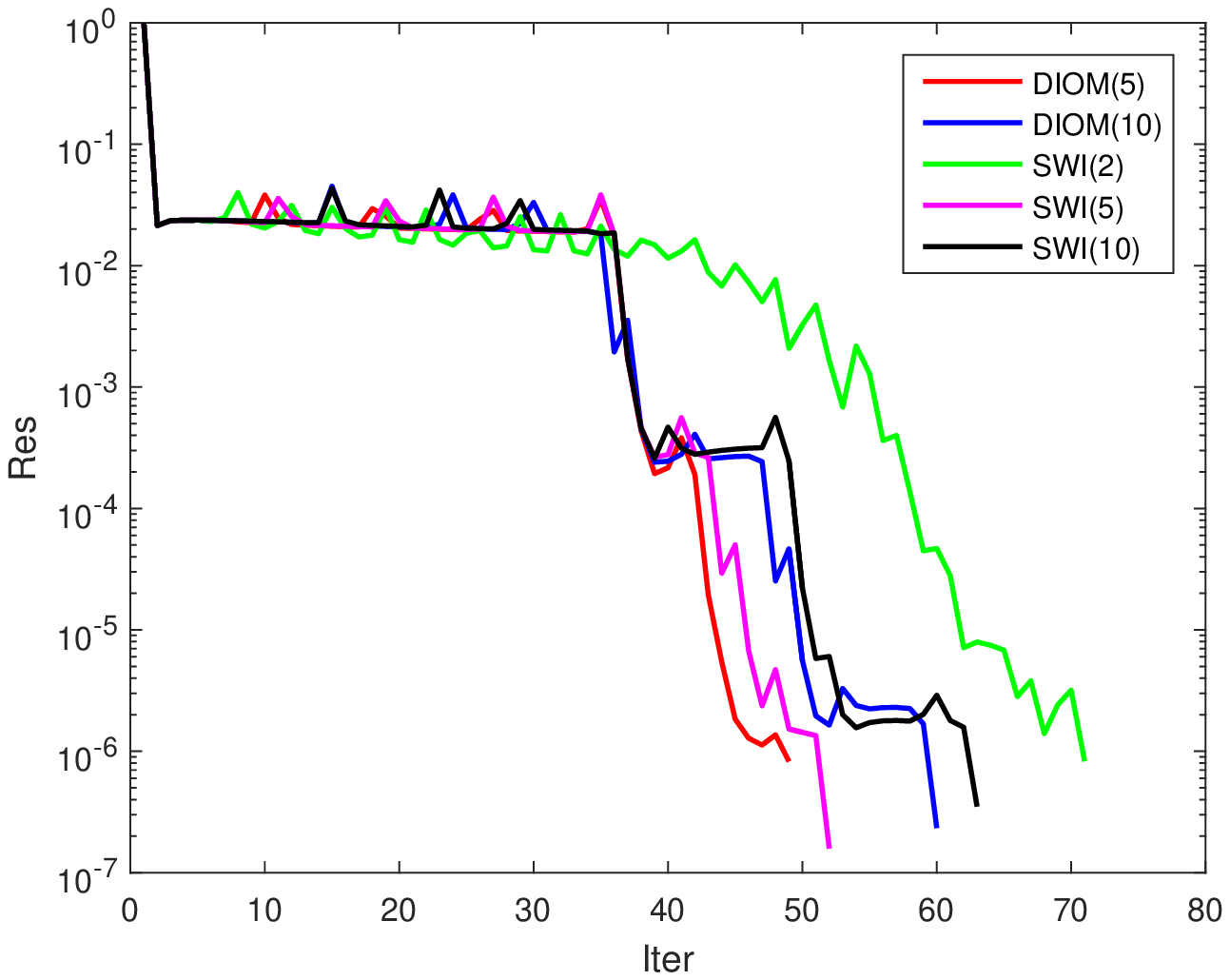}
	}

	\caption{Evolution of the relative residual of the methods tested on \Cref{ex_condiff} with $n=1089$.}
	\label{fig:res}
\end{figure}

In our tests, we set $\epsilon=1/200$ and discretize the convection-diffusion equation using the standard Q1 finite element approximation \citep{Elman2007} on uniform grids with grid parameter $h=1/2^5$, $1/2^6$, $1/2^7$, $1/2^8$, $1/2^9$, $1/2^{10}$. The resulting matrices are unsymmetric positive definite. This discretization was accomplished using IFISS \citep{Elman2007}. We use $m=2$, $5$, and $10$ in SWI, DIOM, and DQGMRES. We report our numerical results in \Cref{fig:res} and \Cref{tab_condiff}.

\begin{table}[H]\scriptsize   
  \centering
  \caption{Numerical results for \Cref{ex_condiff}.}\label{tab_condiff}
    \begin{tabular}{|c|c|c|c|c|c|c|c|c|c|c|c|}
    \hline
$h$&  & $1/2^5$ & $1/2^6$ & $1/2^7$ & $1/2^8$ &$1/2^9$&$1/2^{10}$\\\hline
$n$&  & 1089 &  4225 & 16641 & 66049 & 263169 & 1050625\\\hline

      &Iter  & 43  &  75  & 149  &  298  & 597 & 1193\\
GMRES &CPU   & 0.015 & 0.20  &  1.21  &  14.40  & 347.50 & 36374.39\\
      &Res   & 4.17e-07 & 8.11e-07 & 4.1934e-07 & 7.98e-07  & 7.73e-07 & 8.35e-07\\\hline

      &Iter & 43   & 75  & 149  &  298  & 598 & 1195 \\
FOM   &CPU  & 0.010 & 0.18  & 1.17  & 17.46   & 336.75 & 36299.16 \\
      &Res  & 4.47e-07 & 8.25e-07 & 4.35e-07 & 9.77e-07 & 6.46e-07 & 9.09e-07\\\hline

        &Iter& 49   & 89   & 179  & 363   & 750 & 1544  \\
DIOM(5) &CPU & 0.0039 & 0.049  & 0.42  &   2.81   & 26.86 &  192.49 \\
        &Res & 8.25e-07 & 7.80e-07 & 5.42e-07 & 8.97e-07 & 5.99e-07 & 9.60e-07  \\\hline

        &Iter& 60  &  99 & 166  & 318   & 621  & 1237\\
DIOM(10)&CPU & 0.0065 & 0.062  & 0.75  &  4.01  & 35.74 & 242.13 \\
        &Res & 2.33e-07 & 6.89e-08 & 5.73e-07& 3.87e-07 & 6.74e-07 & 9.26e-07\\\hline

      & Iter & 43  & 75  & 148  & 297   & 594 &  1185\\
SCG & CPU & 0.0036 & 0.043  & 0.47 &  5.25  & 75.31 & 21488.05\\
      & Res & 4.47e-07 & 8.25e-07 & 4.46e-07 & 7.24e-07  &  6.56e-07 &  9.23e-07\\\hline

      &Iter& 71   & 121  & 210  &  417  & 803 & 1598\\
SWI(2)&CPU & 0.0042 & 0.049  & 0.27  &  1.75   & 16.90 & 103.43\\
      &Res & 8.28e-07 & 8.54e-07 & 8.87e-07 & 5.65e-07 & 7.04e-07 & 9.68e-07\\\hline

      &Iter& 52   &  95 &  181 &  360  & 728  & 1435\\
SWI(5)&CPU & 0.0042 & 0.036  & 0.29 & 2.33 & 21.75 &  145.28\\
      &Res & 1.58e-07  & 4.29e-08 &  8.15e-07 & 6.33e-07 &  8.16e-07 & 9.29e-07 \\\hline

       &Iter& 63   &  102   &  178 & 348 & 691 & 1381\\
SWI(10)&CPU & 0.0048 & 0.039 & 0.57 & 3.53 & 26.07 & 205.36 \\
       &Res & 3.50e-07 & 2.97e-07 & 1.87e-07 & 5.50e-07 & 9.58e-07 & 7.70e-07\\\hline

\end{tabular}
\end{table}   

In \Cref{ex_condiff}, BICGSTAB and DIOM(2) were not able to solve problems within $10^4$ iterations when $h\le 1/2^6$ (DIOM(2) also failed when $h=1/2^5$), so we do not report their results in \Cref{tab_condiff}.
DQGMRES failed for $m=2$, $5$, and $10$.
We found that its performance is highly sensitive to the value of $m$. Indeed, while testing other values of \(m\), we observed that the choice of $m=36$ is effective but $m=35$ is not when $h=1/2^5$. The value of $m$ for DQGMRES on this example should not be too much smaller than the number of iterations of GMRES, which is clearly not practical.
\Cref{fig:res} also illustrates that the convergence behaviors of SWI and DIOM are different.

In \Cref{tab_condiff}, the iteration numbers for all tested methods increase in a regular way, each time nearly twice its previous value, but the CPU times rise sharply.  When $h=1/2^{10}$ $(n=1,050,625)$, the CPU times of GMRES, FOM, and SCG are more than $100$ times those of the sliding window versions. When $h \le 1/2^7$, the best performances are by SWI(2). SWI becomes increasingly better as the problem size increases.



\begin{ex}\label{sparse}
We select matrices from the SuiteSparse Matrix Collection \citep{DavisHu-11,sparse} and set $b$ so that the solution is $x_\star = (1, 1, \ldots , 1)$.
\end{ex}

The total number of tested matrices in \Cref{sparse} is $24$, where the matrices arise from applications such as computational fluid dynamics, circuit simulation, directed weighted graphs, optimization, and power networks.
Their name, dimensions and nature are given in \Cref{tab_sparse_size}, where SPD, UPD and UID mean that the matrix is symmetric positive definite, unsymmetric positive definite and unsymmetric indefinite.

\begin{table}[ht]\scriptsize   
  \centering
  \caption{Dimensions and nature of 24 problems in \Cref{sparse}.} \label{tab_sparse_size}

    \begin{tabular}{|r|r|r|r|r|r|r|r|r|}
    \hline

Problem	&	$n$	& Nature  & Problem	&	$n$	& Nature  \\\hline
ACTIVSg10K	&	20000	&	UPD	&	fpga\_dcop\_35	&	1220	&	UPD	\\
ACTIVSg2000	&	4000	&	UPD	&	majorbasis	&	160000	&	UPD	\\
add20	&	2395	&	UPD	&	pde2961	&	2961	&	UPD	\\
add32	&	4960	&	UPD	&	raefsky2	&	3242	&	UID	\\
adder\_dcop\_01	&	1813	&	UPD	&	raefsky4	&	19779	&	SPD	\\
cage12	&	130228	&	UPD	&	raefsky5	&	6316	&	UPD	\\
cage13	&	445315	&	UPD	&	rajat01	&	6833	&	UID	\\
crashbasis	&	160000	&	UID	&	rajat03	&	7602	&	UPD	\\
ex11	&	16614	&	UPD	&	rajat13	&	7598	&	UID	\\
ex18	&	5773	&	UPD	&	rajat16	&	94294	&	UPD	\\
ex19	&	12005	&	UPD	&	rajat27	&	20640	&	UID	\\
ex35	&	19716	&	UPD	&	swang1	&	3169	&	UPD	\\

   \hline
  \end{tabular}
\end{table}   

\begin{figure}[ht]
	\centering
	\subfloat[Iter]{
		\includegraphics[width=.5\linewidth]{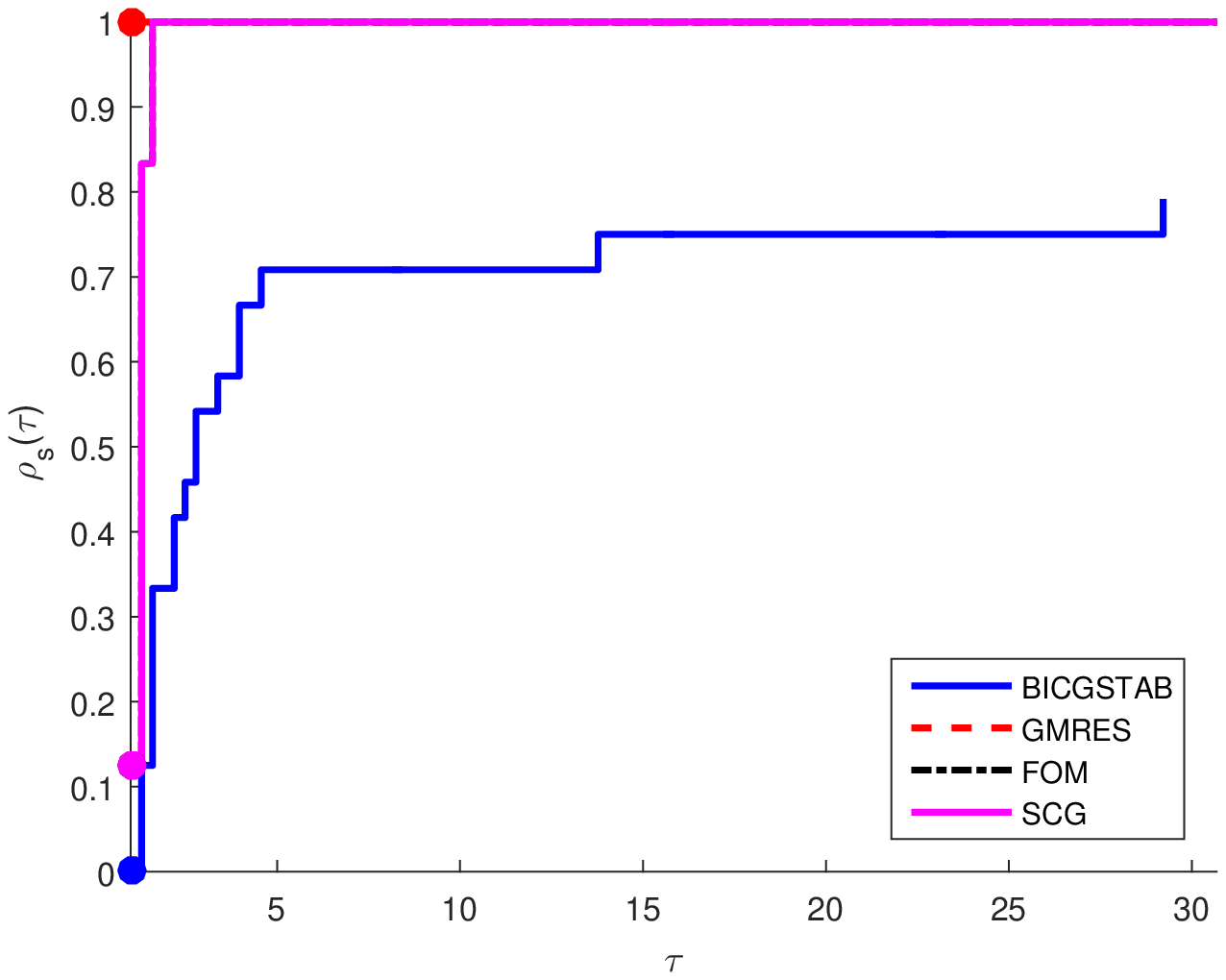}
	}
	\subfloat[CPU]{
		\includegraphics[width=.5\linewidth]{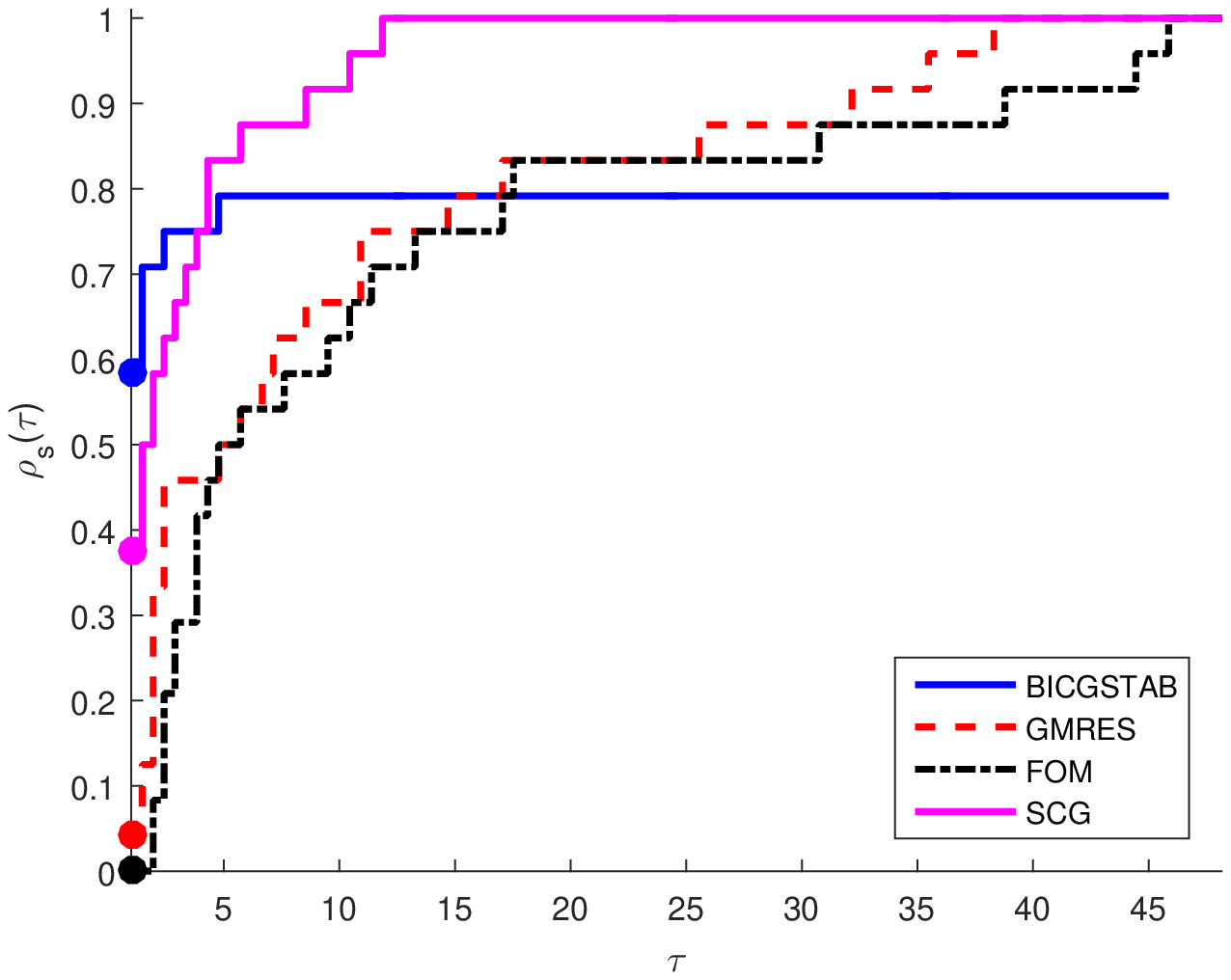}
	}\\
		\subfloat[Iter]{
		\includegraphics[width=.5\linewidth]{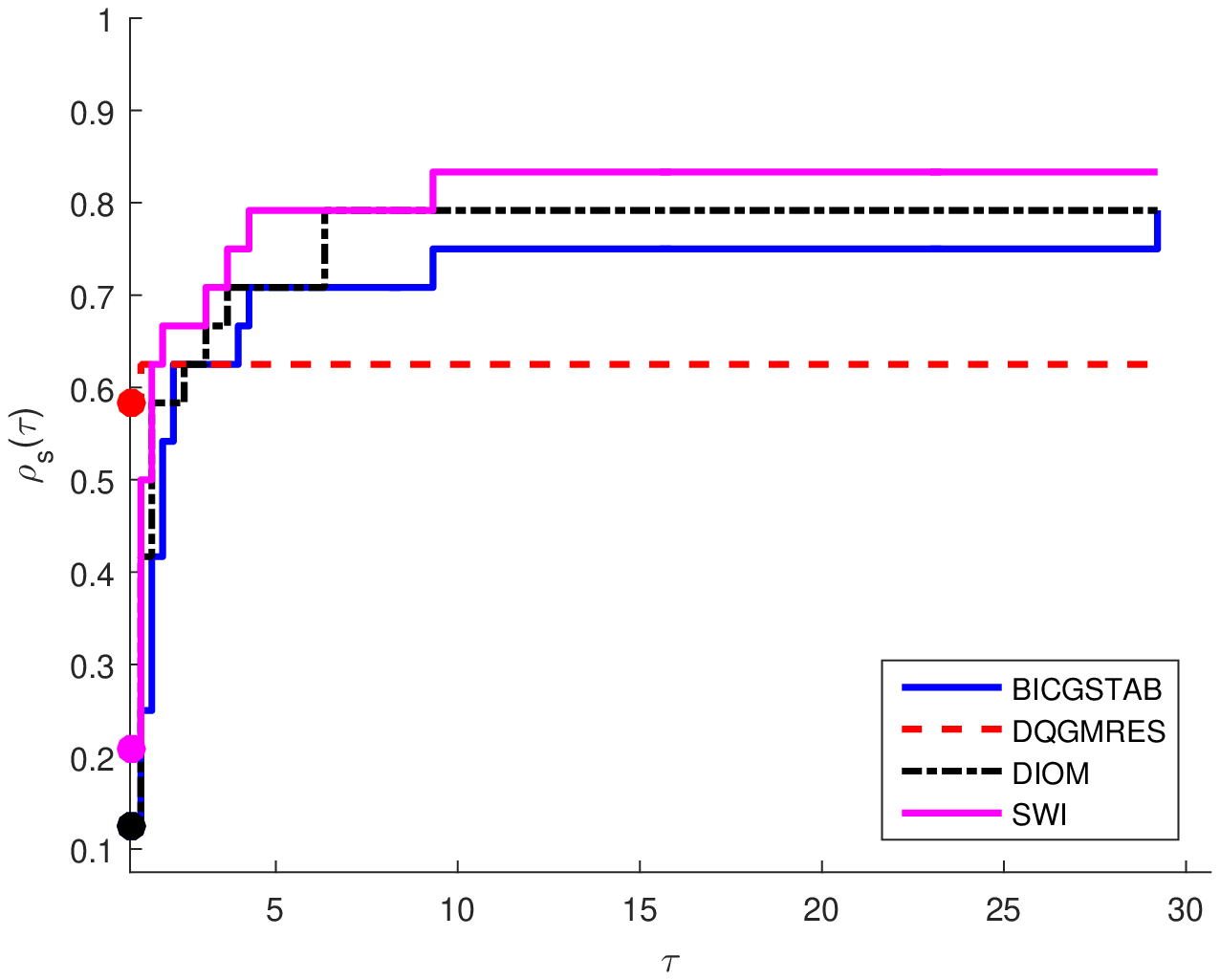}
	}
	\subfloat[CPU]{
		\includegraphics[width=.5\linewidth]{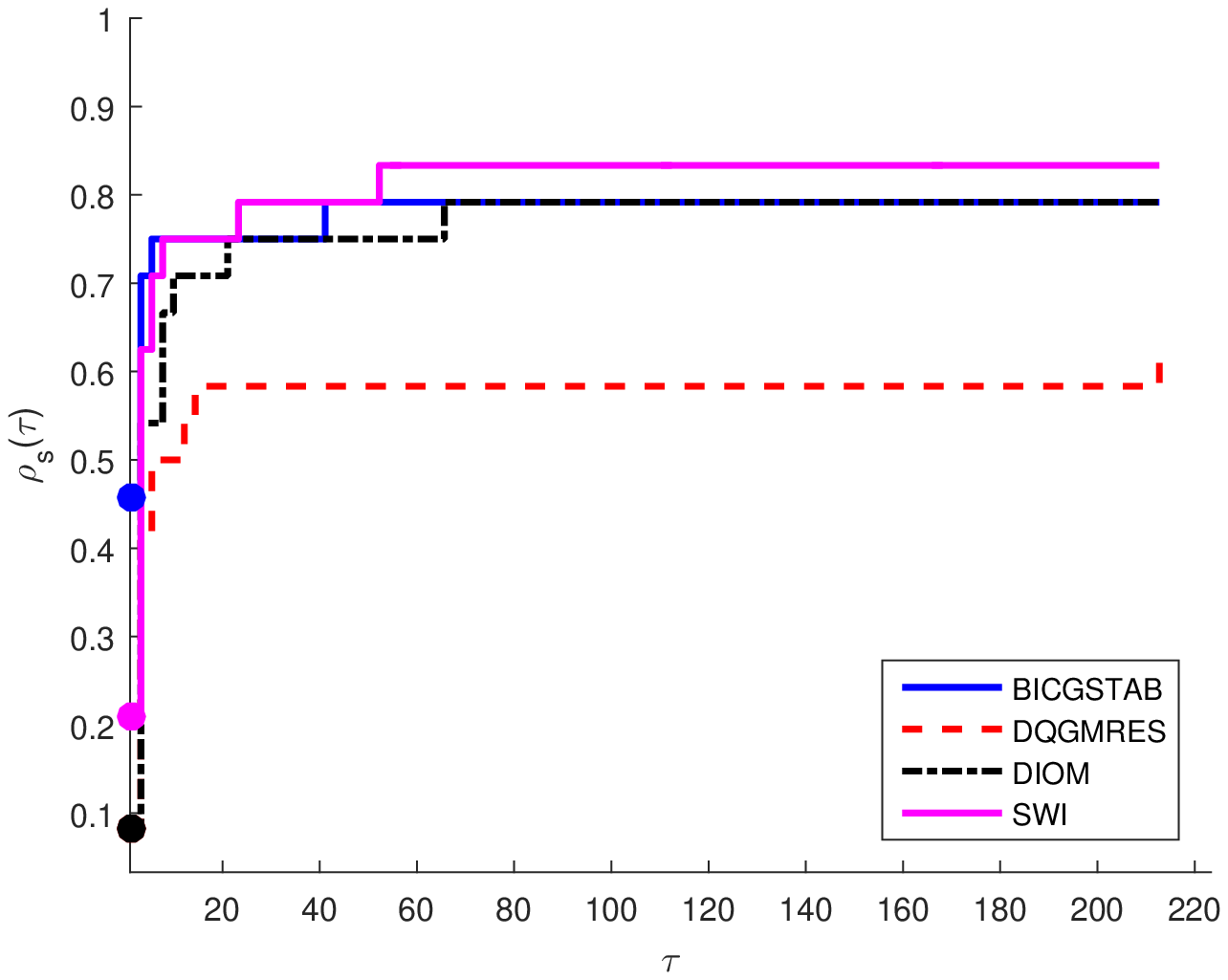}
	}
	\caption{Performance profiles for all tested methods on \Cref{sparse}.
	Limited-memory methods use the value of \(m\) stated in \Cref{tab_sparse_dqgmres,tab_sparse_diom,tab_sparse_swi}.}
	\label{fig:per}
\end{figure}

\begin{figure}[H]
	\centering
	\subfloat[Iter]{
		\includegraphics[width=.5\linewidth]{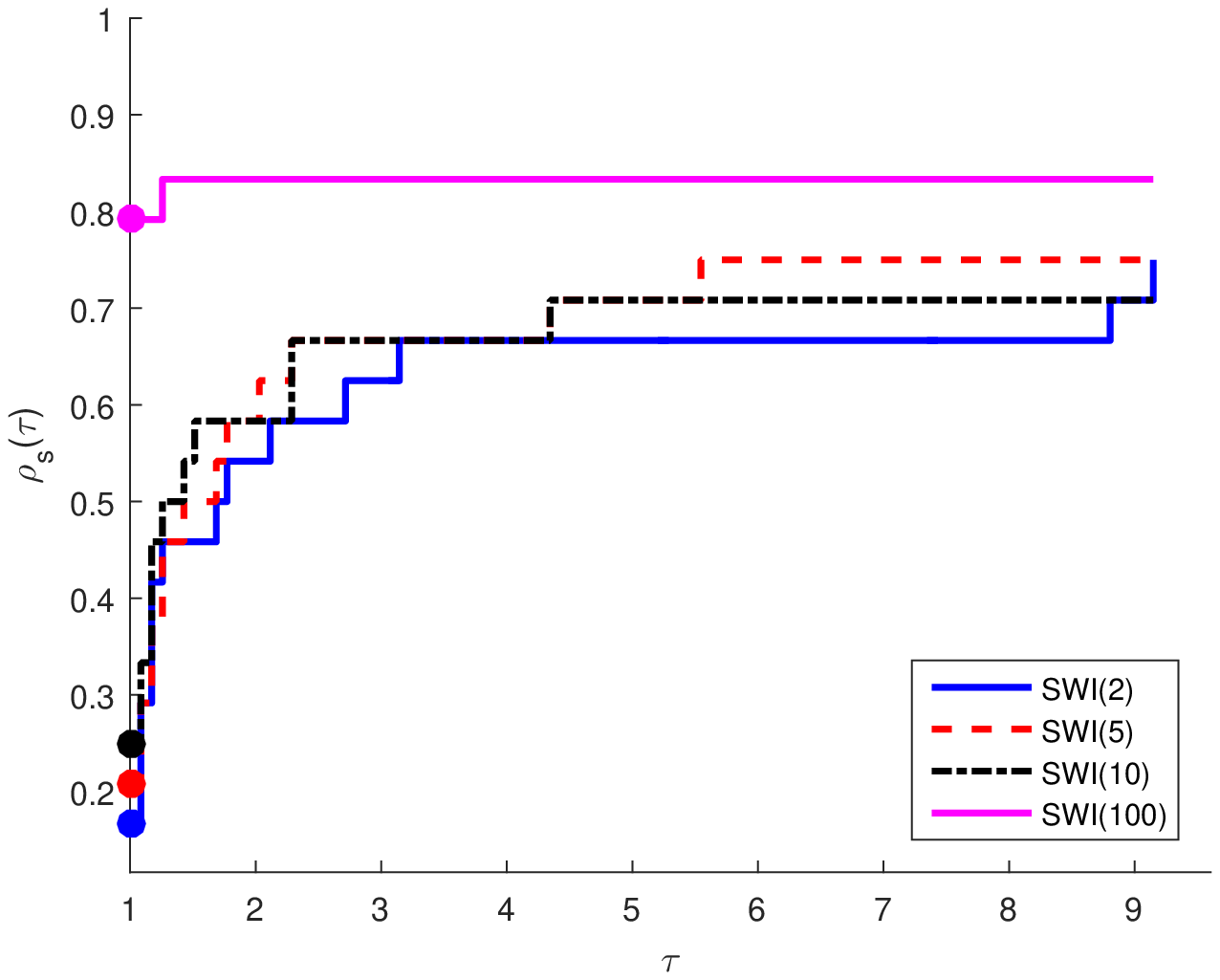}
	}
	\subfloat[CPU]{
		\includegraphics[width=.5\linewidth]{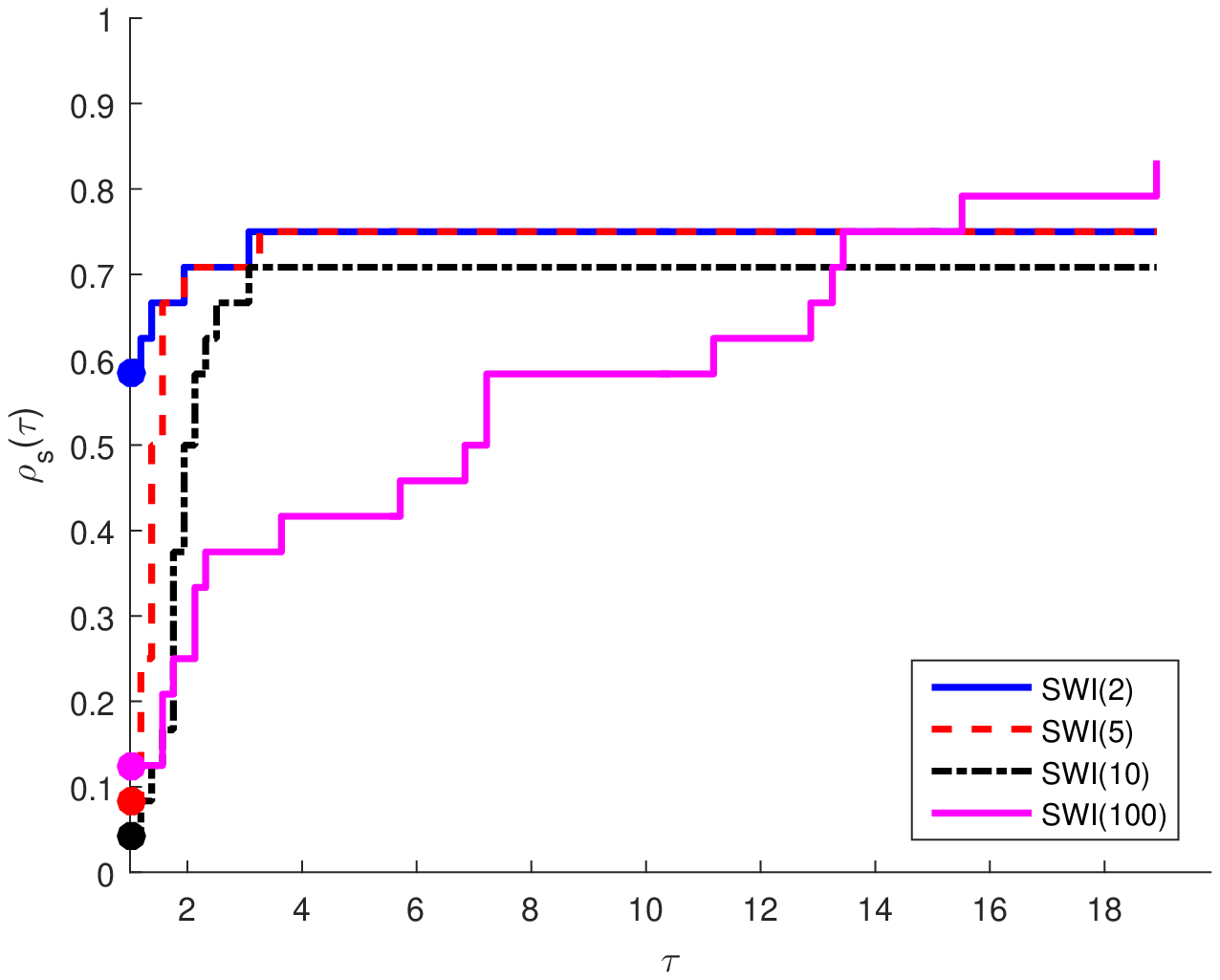}
	}
	\caption{Performance profiles for SWI with different values of $m$ on \Cref{sparse}.}
	\label{fig:per2}
\end{figure}

For DQGMRES, DIOM and SWI, we use $m = 2$, $5$, $10$, and $100$. The best of the four results is presented along with the corresponding value of $m$. The numerical results are reported in \Cref{fig:per,fig:per2} and \Cref{tab_sparse_bicg,tab_sparse_gmres,tab_sparse_dqgmres,tab_sparse_fom,tab_sparse_diom,tab_sparse_scg,tab_sparse_swi}, where ``-" means that the method failed to solve the problem.

GMRES, FOM and SCG successfully solved all the problems but BICGSTAB, DQGMRES, DIOM and SWI failed in 5, 9, 5 and 4 cases, respectively. In terms of the CPU time, BICGSTAB, SCG and SWI perform best in 10, 7, and 4 cases, respectively.
Compared to SCG, SWI requires less CPU time in 14 cases and the improvements are significant. Compared to DQGMRES and DIOM, SWI requires the least CPU time in 16 cases, while DQGMRES and DIOM require the least in only 2 cases. Hence, SWI is the most successful of the sliding window implementations.

In \Cref{fig:per}, performance profiles%
\footnote{The performance profile $\varrho_s(\tau)$ is a distribution function for the performance ratio $r_{p,s}$, with
$$r_{p,s} = t_{p,s} / \min\{ t_{p,s}:\,s\in S\}
\qquad \mbox{and} \qquad
\varrho_s(\tau) = \left| \left\{p\in P:\,r_{p,s}\le \tau\right\} \right| / |P|,
$$
where $S$ is the set of solvers, $P$ is the set of problems, \(|\cdot|\) indicates cardinality, and $t_{p,s}$ is the Iter/CPU required to solve problem $p$ with solver $s$.}
\citep{profiles} indicate that SCG and SWI are more robust than BICGSTAB, and also more efficient than other tested methods;  the reduction in CPU time for SWI was often substantial.
To see the role of $m$ in SWI's performance, we also plot performance profiles for SWI with different $m$. From \Cref{fig:per2} it is apparent that larger $m$ leads to fewer iterations for SWI but more CPU time. Hence, the choice of $m$ to balance iterations and CPU time is crucial for the performance of SWI.

\begin{table}[ht]\scriptsize   
\setlength\tabcolsep{5.2pt}
  \centering
  \caption{Numerical results for BICGSTAB on \Cref{sparse}.}\label{tab_sparse_bicg}
    \begin{tabular}{|r|r|r|r|r|r|r|r|}
    \hline

Problem	&	Iter	&	CPU	&	Res	&	Problem	&	Iter	&	CPU	&	Res	\\ \hline
ACTIVSg10K	&	-	&	-	&	-	&	fpga\_dcop\_35	&	-	&	- 	&	-	\\
ACTIVSg2000	&	2608	&	0.42 	&	3.04E-07	&	majorbasis	&	110	&	0.71 	&	8.20E-07	\\
add20	&	376	&	0.053 	&	8.59E-07	&	pde2961	&	267	&	0.021 	&	9.12E-07	\\
add32	&	72	&	0.0086 	&	8.44E-07	&	raefsky2	&	636	&	0.32 	&	6.92E-07	\\
adder\_dcop\_01	&	1607	&	0.072 	&	9.96E-07	&	raefsky4	&	29	&	0.086 	&	7.01E-07	\\
cage12	&	20	&	0.11 	&	5.71E-07	&	raefsky5	&	129	&	0.036 	&	7.49E-07	\\
cage13	&	20	&	0.40 	&	4.76E-07	&	rajat01	&	-	&	-	&	-	\\
crashbasis	&	244	&	1.57 	&	9.74E-07	&	rajat03	&	2338	&	0.56 	&	4.08E-07	\\
ex11	&	1572	&	3.56 	&	8.94E-07	&	rajat13	&	86	&	0.014 	&	9.52E-07	\\
ex18	&	1397	&	0.25 	&	9.89E-07	&	rajat16	&	-	&	-	&	-	\\
ex19	&	4050	&	2.88 	&	3.79E-07	&	rajat27	&	-	&	-	&	-	\\
ex35	&	1438	&	1.17 	&	8.15E-07	&	swang1	&	22	&	0.0017 	&	6.14E-07	\\

   \hline
  \end{tabular}
\end{table}   

\begin{table}[ht]\scriptsize   
\setlength\tabcolsep{5.2pt}
  \centering
  \caption{Numerical results for GMRES on \Cref{sparse}.}\label{tab_sparse_gmres}
    \begin{tabular}{|r|r|r|r|r|r|r|r|}
    \hline
Problem	&	Iter	&	CPU	&	Res	&	Problem	&	Iter	&	CPU	&	Res	\\ \hline
ACTIVSg10K	&	6130	&	2690.28 	&	9.98E-07	&	fpga\_dcop\_35	&	214	&	0.29 	&	9.37E-07	\\
ACTIVSg2000	&	779	&	16.03 	&	9.71E-07	&	majorbasis	&	97	&	7.49 	&	9.81E-07	\\
add20	&	195	&	0.35 	&	9.92E-07	&	pde2961	&	189	&	0.67 	&	9.11E-07	\\
add32	&	57	&	0.091 	&	9.42E-07	&	raefsky2	&	332	&	2.71 	&	8.61E-07	\\
adder\_dcop\_01	&	55	&	0.023 	&	9.39E-07	&	raefsky4	&	22	&	0.085 	&	4.73E-07	\\
cage12	&	14	&	0.16 	&	9.57E-07	&	raefsky5	&	34	&	0.048 	&	8.49E-07	\\
cage13	&	15	&	0.65 	&	8.33E-07	&	rajat01	&	1894	&	132.19 	&	9.99E-07	\\
crashbasis	&	175	&	22.80 	&	9.93E-07	&	rajat03	&	172	&	1.32 	&	8.98E-07	\\
ex11	&	571	&	24.24 	&	9.98E-07	&	rajat13	&	22	&	0.022 	&	9.06E-07	\\
ex18	&	505	&	8.75 	&	9.99E-07	&	rajat16	&	1116	&	355.88 	&	9.97E-07	\\
ex19	&	936	&	48.60 	&	9.31E-07	&	rajat27	&	595	&	30.20 	&	9.89E-07	\\
ex35	&	607	&	29.71 	&	9.99E-07	&	swang1	&	19	&	0.0076 	&	7.15E-07	\\

   \hline
  \end{tabular}
\end{table}   

\begin{table}[H]\scriptsize   
\setlength\tabcolsep{5.2pt}
  \centering
  \caption{Numerical results for FOM on \Cref{sparse}.}\label{tab_sparse_fom}
    \begin{tabular}{|r|r|r|r|r|r|r|r|}
    \hline

Problem	&	Iter	&	CPU	&	Res	&	Problem	&	Iter	&	CPU	&	Res	\\ \hline
ACTIVSg10K	&	6376	&	2762.76 	&	9.48E-07	&	fpga\_dcop\_35	&	276	&	0.42 	&	5.26E-07	\\
ACTIVSg2000	&	814	&	18.55 	&	9.75E-07	&	majorbasis	&	109	&	9.36 	&	8.54E-07	\\
add20	&	214	&	0.40 	&	8.59E-07	&	pde2961	&	192	&	0.63 	&	8.36E-07	\\
add32	&	59	&	0.095 	&	6.06E-07	&	raefsky2	&	334	&	2.93 	&	9.72E-07	\\
adder\_dcop\_01	&	80	&	0.055 	&	9.41E-07	&	raefsky4	&	22	&	0.16 	&	5.06E-07	\\
cage12	&	15	&	0.27 	&	5.12E-07	&	raefsky5	&	36	&	0.077 	&	7.37E-07	\\
cage13	&	15	&	0.84 	&	9.47E-07	&	rajat01	&	2296	&	185.61 	&	5.60E-07	\\
crashbasis	&	196	&	27.43 	&	9.89E-07	&	rajat03	&	173	&	1.28 	&	7.55E-07	\\
ex11	&	694	&	35.94 	&	9.96E-07	&	rajat13	&	24	&	0.025 	&	7.72E-07	\\
ex18	&	541	&	9.64 	&	9.62E-07	&	rajat16	&	1453	&	580.32 	&	9.89E-07	\\
ex19	&	947	&	47.85 	&	9.41E-07	&	rajat27	&	909	&	68.95 	&	9.45E-07	\\
ex35	&	847	&	53.71 	&	9.35E-07	&	swang1	&	19	&	0.0075 	&	8.26E-07	\\

   \hline
  \end{tabular}
\end{table}   

\begin{table}[H]\scriptsize   
\setlength\tabcolsep{5.2pt}
  \centering
  \caption{Numerical results for SCG on \Cref{sparse}.}\label{tab_sparse_scg}
    \begin{tabular}{|r|r|r|r|r|r|r|r|}
    \hline

Problem	&	Iter	&	CPU	&	Res	&	Problem	&	Iter	&	CPU	&	Res	\\ \hline
ACTIVSg10K	&	6376	&	1266.00 	&	9.50E-07	&	fpga\_dcop\_35	&	276	&	0.13 	&	5.26E-07	\\
ACTIVSg2000	&	814	&	4.34 	&	9.75E-07	&	majorbasis	&	109	&	2.50 	&	8.54E-07	\\
add20	&	214	&	0.084 	&	8.59E-07	&	pde2961	&	192	&	0.12 	&	8.36E-07	\\
add32	&	59	&	0.020 	&	6.06E-07	&	raefsky2	&	334	&	0.57 	&	9.72E-07	\\
adder\_dcop\_01	&	80	&	0.016 	&	9.41E-07	&	raefsky4	&	22	&	0.094 	&	5.06E-07	\\
cage12	&	15	&	0.14 	&	5.12E-07	&	raefsky5	&	36	&	0.027 	&	7.37E-07	\\
cage13	&	15	&	0.45 	&	9.47E-07	&	rajat01	&	2294	&	70.42 	&	9.75E-07	\\
crashbasis	&	196	&	6.14 	&	9.89E-07	&	rajat03	&	173	&	0.24 	&	7.60E-07	\\
ex11	&	694	&	8.97 	&	9.96E-07	&	rajat13	&	24	&	0.013 	&	7.72E-07	\\
ex18	&	541	&	2.08 	&	9.62E-07	&	rajat16	&	1453	&	158.62 	&	9.98E-07	\\
ex19	&	947	&	12.28 	&	9.41E-07	&	rajat27	&	912	&	16.63 	&	7.96E-07	\\
ex35	&	847	&	13.47 	&	9.33E-07	&	swang1	&	19	&	0.0056 	&	8.26E-07	\\

   \hline
  \end{tabular}
\end{table}   

\begin{table}[ht]\scriptsize   
\setlength\tabcolsep{3.5pt}
  \centering
  \caption{Numerical results for DQGMRES on \Cref{sparse}.}\label{tab_sparse_dqgmres}
    \begin{tabular}{|r|r|r|r|r|r|r|r|r|r|r|}
    \hline

Problem	&	Iter	&	CPU	&	Res	& $m$ &	Problem	&	Iter	&	CPU	&	Res	& $m$ \\ \hline
ACTIVSg10K	&	-	&	-	&	-	&	-	&	fpga\_dcop\_35	&	97	&	7.35	&	9.81E-07	&	100	\\
ACTIVSg2000	&	-	&	-	&	-	&	-	&	majorbasis	&	-	&	-	&	-	&	-	\\
add20	&	206	&	0.46 	&	9.98E-07	&	100	&	pde2961	&	-	&	-	&	-	&	-	\\
add32	&	57	&	0.094 	&	9.42E-07	&	100	&	raefsky2	&	-	&	-	&	-	&	-	\\
adder\_dcop\_01	&	55	&	0.022 	&	9.39E-07	&	100	&	raefsky4	&	22	&	0.093 	&	4.73E-07	&	5	\\
cage12	&	15	&	0.25 	&	5.65E-07	&	5	&	raefsky5	&	34	&	0.045 	&	8.49E-07	&	100	\\
cage13	&	16	&	0.79 	&	7.25E-07	&	5	&	rajat01	&	-	&	-	&	-	&	-	\\
crashbasis	&	-	&	-	&	-	&	-	&	rajat03	&	257	&	0.17 	&	9.86E-07	&	5	\\
ex11	&	893	&	4.87 	&	9.87E-07	&	5	&	rajat13	&	45	&	0.071 	&	8.34E-07	&	10	\\
ex18	&	798	&	0.53 	&	9.90E-07	&	5	&	rajat16	&	-	&	-	&	-	&	-	\\
ex19	&	1897	&	3.54 	&	9.78E-07	&	5	&	rajat27	&	-	&	-	&	-	&	-	\\
ex35	&	347	&	0.029 	&	9.83E-07	&	5	&	swang1	&	19	&	0.0080 	&	7.44E-07	&	10	\\

   \hline
  \end{tabular}
\end{table}   

\begin{table}[ht]\scriptsize   
\setlength\tabcolsep{3.5pt}
  \centering
  \caption{Numerical results for DIOM on \Cref{sparse}.}\label{tab_sparse_diom}
    \begin{tabular}{|r|r|r|r|r|r|r|r|r|r|r|}
    \hline

Problem	&	Iter	&	CPU	&	Res	& $m$ &	Problem	&	Iter	&	CPU	&	Res	& $m$ \\ \hline
ACTIVSg10K	&	-	&	-	&	-	&	-	&	fpga\_dcop\_35	&	587	&	0.035 	&	6.07E-07	&	2	\\
ACTIVSg2000	&	-	&	-	&	-	&	-	&	majorbasis	&	268	&	3.98 	&	9.83E-07	&	2	\\
add20	&	224	&	0.038 	&	9.96E-07	&	5	&	pde2961	&	384	&	0.18 	&	9.58E-07	&	10	\\
add32	&	59	&	0.015 	&	6.26E-07	&	2	&	raefsky2	&	3883	&	6.05 	&	9.63E-07	&	10	\\
adder\_dcop\_01	&	80	&	0.049 	&	9.41E-07	&	100	&	raefsky4	&	22	&	0.091 	&	5.06E-07	&	2	\\
cage12	&	15	&	0.17 	&	5.86E-07	&	2	&	raefsky5	&	54	&	0.068 	&	8.50E-07	&	10	\\
cage13	&	16	&	0.52 	&	4.63E-07	&	2	&	rajat01	&	-	&	-	&	-	&	-	\\
crashbasis	&	721	&	10.91 	&	9.96E-07	&	2	&	rajat03	&	281	&	0.12 	&	6.78E-07	&	2	\\
ex11	&	1237	&	6.04 	&	9.58E-07	&	2	&	rajat13	&	55	&	0.10 	&	9.54E-07	&	10	\\
ex18	&	943	&	0.45 	&	9.80E-07	&	2	&	rajat16	&	-	&	-	&	-	&	-	\\
ex19	&	2181	&	2.97 	&	8.87E-07	&	2	&	rajat27	&	-	&	-	&	-	&	-	\\
ex35	&	1235	&	1.90 	&	8.03E-07	&	2	&	swang1	&	19	&	0.0049 	&	9.46E-07	&	5	\\

   \hline
  \end{tabular}
\end{table}   

\begin{table}[H]\scriptsize   
\setlength\tabcolsep{3.5pt}
  \centering
  \caption{Numerical results for SWI on \Cref{sparse}.}\label{tab_sparse_swi}
    \begin{tabular}{|r|r|r|r|r|r|r|r|r|r|r|}
    \hline

Problem	&	Iter	&	CPU	&	Res	& $m$ &	Problem	&	Iter	&	CPU	&	Res	& $m$ \\ \hline
ACTIVSg10K	&	-	&	-	&	-	&	-	&	fpga\_dcop\_35	&	903	&	0.053 	&	4.54E-07	&	2	\\
ACTIVSg2000	&	1764	&	9.35 	&	9.94E-07	&	100	&	majorbasis	&	130	&	1.44 	&	9.78E-07	&	2	\\
add20	&	234	&	0.033 	&	9.46E-07	&	2	&	pde2961	&	255	&	0.076 	&	9.21E-07	&	10	\\
add32	&	59	&	0.014 	&	6.23E-07	&	2	&	raefsky2	&	2683	&	1.95 	&	9.95E-07	&	5	\\
adder\_dcop\_01	&	80	&	0.037 	&	9.41E-07	&	100	&	raefsky4	&	22	&	0.059 	&	5.06E-07	&	2	\\
cage12	&	15	&	0.13 	&	5.28E-07	&	2	&	raefsky5	&	95	&	0.046 	&	8.29E-07	&	2	\\
cage13	&	15	&	0.40 	&	9.60E-07	&	2	&	rajat01	&	-	&	-	&	-	&	-	\\
crashbasis	&	428	&	7.52 	&	9.72E-07	&	5	&	rajat03	&	341	&	0.12 	&	7.50E-07	&	2	\\
ex11	&	1235	&	3.38 	&	9.05E-07	&	2	&	rajat13	&	53	&	0.038 	&	9.86E-07	&	10	\\
ex18	&	908	&	0.36 	&	9.64E-07	&	2	&	rajat16	&	-	&	-	&	-	&	-	\\
ex19	&	2188	&	2.22 	&	8.52E-07	&	2	&	rajat27	&	-	&	-	&	-	&	-	\\
ex35	&	1242	&	1.48 	&	9.36E-07	&	2	&	swang1	&	20	&	0.0038 	&	9.31E-07	&	2	\\

   \hline
  \end{tabular}
\end{table}   



\section{Conclusions and future work}\label{sec:con}

We introduced the semi-conjugate gradient method (SCG) and its sliding window implementation (SWI) to solve unsymmetric positive definite linear systems. Both theoretical and numerical studies of SCG and SWI were conducted. SCG is theoretically equivalent to FOM, but a counter-example illustrates that their sliding window implementations differ. The numerical results presented are highly encouraging, though the performance of SWI naturally depends on the window width $m$.

Future work should aim to develop efficient algorithms for adaptive selection of the window width $m$. A possibly feasible approach is changing the value of $m$ dynamically. Another way to improve the performance of SCG and SWI is to incorporate preconditioning into them and to develop practical and effective preconditioners.
It is also interesting and challenging to extend SCG and SWI to solve nonlinear problems as has been done for CG
\cite{Oleg2018}.

\section*{Acknowledgments}

We would like to thank our colleague and friend, Dr Oleg Burdakov, for his devotion to research and for his everlasting sense of humor.
In particular, we wish to express our gratitude to him for fundamental contributions that initiated this work, and for many constructive suggestions on our early Matlab implementation of SCG and SWI.
At the end of 2017, he independently constructed SCG and SWI to solve quasi-definite linear systems. Subsequently, he extended them to general unsymmetric positive definite linear systems.  The question of choosing an ideal $m$ remains for future research, as it does for GMRES($m$) and DQGMRES($m$).



\bibliographystyle{abbrvnat}
\bibliography{references}

\end{document}